\documentclass{amsart}
\usepackage{amssymb,bbm,calligra,mathrsfs}
\usepackage[all]{xy}

\def\calli#1{\textup{\!\textcalligra{#1}}}
\def\tl{{\underline\circledcirc}}
\def\op{^{\mathrm{op}}}

\def\1{^{-1}}

\def\xto#1{\xrightarrow[]{#1}}

\DeclareMathOperator\id{id} \DeclareMathOperator\Ob{Ob}
 \DeclareMathOperator\Aut{Aut}
 
\DeclareMathOperator\Ch{Ch}

\DeclareMathOperator\Tracks{Trext} 
\DeclareMathOperator\Xext{Xext} \DeclareMathOperator\pair{pair}

\def\ubx{{\underline\boxtimes}}

\let\x\times
\let\ho\simeq

\let\ot\leftarrow

\let\d\partial
\let\then\Rightarrow

\let\ge\geqslant

\def\1{^{-1}}

\def\tP{{\tilde{P}}}

\def\f{{\mathbf f}}

\def\xto#1{\xrightarrow[]{#1}}

\def\brk#1{\langle#1\rangle}

\def\cref#1#2#3{\left(#2\right.\left|\ #3\right)_{#1}}

\newtheorem{Pro}{Proposition}[subsection]
\newtheorem{Le}[Pro]{Lemma}
\newtheorem{The}[Pro]{Theorem}
\newtheorem{Co}[Pro]{Corollary}
\theoremstyle{definition}
\newtheorem{De}[Pro]{Definition}
\theoremstyle{remark}
\newtheorem{Exm}[Pro]{Example}
\newtheorem{Rem}[Pro]{Remark}

\let\bb\beta

\def \Im{\mathop{\sf Im}\nolimits}
\def \Ker{\mathop{\sf Ker}\nolimits}
\def \Hom{\mathop{\sf Hom}\nolimits}
\def \Ext{\mathop{\sf Ext}\nolimits}

\newcommand{\Z}{{\mathbb Z}}

\let\aa\alpha

\def\ta{{\cal T}}
\let\ho\simeq

\def\c{{\cal C}}
\def\ka{{\cal K}}
\def\Nil{{\sf Nil}}
\def\QPM{{\sf{qpm}}}

\let\x\times
\def\modf{\mathbf{mod}\textrm{-}}
\let\tp\otimes
\def\La{\Lambda}
\def\Mod{\mathbf{Mod}\textrm{-}}

\let\t\otimes

\def\Tor{\mathop{\sf Tor}\nolimits}

\def\Cok{\mathop{\sf Coker}\nolimits}
\def\H{\mathop{\sf H}\nolimits}
\def\Sh{\mathop{\sf SH}\nolimits}
\def\HML{\mathop{\sf HML}\nolimits}

\let\cal\mathscr

\let\ge\geqslant
\let\le\leqslant

\let\x\times
\let\ox\otimes

\let\ot\leftarrow
\let\then\Rightarrow

\let\then\Rightarrow

\def\xto#1{\xrightarrow[]{#1}}

\def\ab{_{\mathrm{ab}}}

\def\brk#1{\left\langle#1\right\rangle}

\def\hog#1{[\![#1]\!]}

\def\id{{\mathrm{Id}}}

\def\op{^{\mathrm{op}}}
\def\1{^{-1}}

\def\Tracks{\sf Tracks}

\def\A{{\mathbb A}}
\def\B{{\mathbb B}}

\def\D{{\mathbb D}}

\def\C{{\mathbb C}}

\def\N{{\mathbb N}}

\def\Z{{\mathbb Z}}

\def\QR{{\sf QR}}
\def\SG{{\sf SG}}
\def\SR{{\sf SR}}

\def\CSG{{\sf cross(SG)}}

\def\ad{^{\mathrm{ad}}}

\def\c{{\mathscr C}}

\def\f{{\mathscr F}}

\def\V{{\mathbf V}}

\def\F{{\sf F}}
\def\Ab{{\sf Ab}}
\def\Track{{\Tracks}}

\def\CS{{\sf Csr}}
\def\Gr{{\sf Gr}}
\def\Sets{{\sf Sets}}
\def\Gpd{{\sf Gpd}}
\def\Cat{{\sf Cat}}

\def\Nil{{\sf Nil}}

\def\Id{{\sf Id}}
\def\modr{\mathbf{mod}\textrm{-}}
\def\Modr{\mathbf{Mod}\textrm{-}}

\numberwithin{equation}{subsection}

\def\xext{{\calli{Xext}\;}}

\def\Znil{\Z_{\mathrm{nil}}}

\begin{document}

\title{Third Mac Lane cohomology}
\author{H.-J. Baues}
\address{
Max-Planck-Institut f\"ur Mathematik\\
Vivatsgasse 7\\
Bonn 53111\\
Germany} \email{baues@mpim-bonn.mpg.de}

\author{M. Jibladze}
\address{
A. Razmadze Mathematical Institute\\
M. Alexidze st. 1\\
Tbilisi 0193\\
Georgia} \email{jib@rmi.acnet.ge}

\author{T. Pirashvili}
\address{
A. Razmadze Mathematical Institute\\
M. Alexidze st. 1\\
Tbilisi 0193\\
Georgia} \email{pira@rmi.acnet.ge}

\begin{abstract}MacLane cohomology is an algebraic version of the topological
Hochschild cohomology.
Based on the computation of the third author (see Appendix below) we
obtain  an interpretation of the third Mac~Lane cohomology
of rings
using certain kind of crossed extensions of rings in the quadratic
world. Actually we obtain two such interpretations corresponding to
the two monoidal structures on the category of square groups.
\end{abstract}

\maketitle

\section{Introduction}\label{intro}

Let $R$ be a ring and $M$ a bimodule over $R$. Then there are three
essential cohomology theories associated to the pair $(R,M)$ due to
Hochschild, Shukla, and Mac~Lane, see \cite{hh, sh, macl}. These theories are connected by
natural maps (\cite{shukla})
$$
\H^n(R;M)\to\Sh^n(R;M)\xto{\tau^n}\HML^n(R;M).
$$
It is known that Mac~Lane cohomology coincides with topological
Hochschild cohomology (\cite{fw}) and coincides also with
Baues-Wirsching cohomology of the category $\modr R$ of finitely
generated free $R$-modules (\cite{JP}). We study the cohomologies in
dimension $n=3$. In this case the elements in $\H^3(R;M)$ are
represented by split crossed extensions and elements in $\Sh^3(R;M)$
are represented by all crossed extensions of $R$ by $M$ in the
monoidal category $(\Ab,\ox)$ of abelian groups. See \cite{jll},
\cite{baumin} and \cite{shukla}. Here a crossed extension of $R$ by
$M$ is an exact sequence in $\Ab$,
$$
0\to M\xto\iota C_1\xto\d C_0\xto q R\to0,
$$
where $C_0$ is a ring, $C_1$ is a bimodule over it, $q$ is a ring
homomorphism and $\iota$ and $\d$ are $C_0$-biequivariant maps
satisfying $b\d(c)=\d(b)c$ for $b,c\in C_1$.

A similar result for group cohomology in dimension 3 is due to
Mac~Lane and J.~H.~C.~Whitehead \cite{MW}.

The main goal of this paper is the construction of appropriate
crossed extensions of $R$ by $M$ which represent classes in the
third Mac~Lane cohomology $\HML^3(R;M)$.

To this end we recall that for any small category $\c$ with a
natural system $D$ on it, the Baues-Wirsching cohomology group
$H^3(\c;D)$ can be represented by linear track extensions of $\c$ by
$D$, see \cite{P1,P2,B}. Hence by the isomorphism
$$
\HML^3(R;M)\cong H^3(\modr R;\Hom_R(-,-\ox_RM))
$$
elements of $\HML^3(R;M)$ are represented by linear track extensions
of $\modr R$. Such a description, however, is available for any
category $\c$ and does not restrict to the specific nature of
Mac~Lane cohomology of a ring.

In order to find specific crossed extensions for $\HML^3(R;M)$ we
have to proceed from linear algebra to quadratic algebra. Here
``linear algebra'' is the algebra of rings and modules. A ring is a monoid
in the monoidal category $(\Ab,\ox)$ and a module is an object in
$\Ab$ together with an action of such a monoid.

In ``quadratic algebra'' abelian groups are replaced by square
groups. In fact, if one considers endofunctors of the category of
groups which preserve filtered colimits and reflexive coequalizers,
then abelian groups can be identified with linear endofunctors and
square groups can be identified with quadratic endofunctors
(\cite{square}). The category $\SG$ of square groups contains the
category $\Ab$ as a full subcategory since a linear endofunctor is
also quadratic. Composition of functors leads to monoidal structures
$\t$ and $\square$ in such a way that $(\Ab,\t)$ is a monoidal
subcategory of $(\SG,\square)$. There is also another monoidal structure
$\tl$ on $\SG$ such that the identity of $\SG$ is a lax monoidal
functor $(\SG,\tl)\to(\SG,\square)$. Here $\tl$ is symmetric while $\square$
is highly nonsymmetric. Compare \cite{bjp}.

Crossed extensions in the monoidal categories $(\SG,\square)$ or
$(\SG,\tl)$ are defined similarly to the case $(\Ab,\t)$ above, see
section \ref{crex}. As a main result we prove in this paper the
following theorem, compare the more detailed version \ref{xext}
below.

\subsection{Main theorem}

\begin{The}\label{mainth}
Elements in the third Mac~Lane cohomology group $\HML^3(R;M)$ are in
1-1 correspondence with equivalence classes of linearly generated
crossed extensions of $R$ by $M$ in the monoidal category
$(\SG,\tl)$, or in the monoidal category $(\SG,\square)$.
\end{The}

Such an interpretation of the group $\HML^3$ was missing for many
years; in terms of obstruction theory the problem first arose in the
classical paper of Mac~Lane \cite{maclo}. The theorem is based on
the quadratic theory developed in \cite{square,bjp} and emphasizes
importance of the quadratic algebra of square groups. A crucial step
in the proof of the theorem relies on the vanishing result achieved
by the third named author in the Appendix.

In dimension three, the map $\tau^3$ fits in the exact sequence (see
\cite{JPmoambe}, \cite{shukla})
$$
0\to\Sh^3(R;M)\xto{\tau^3}\HML^3(R;M)\xto\nu\H^0(R;{}_2M)\to
\Sh^4(R;M)\xto{\tau^4}\HML^4(R;M)
$$
where $_2M=\{m\in M\mid 2m=0\}$.

As an application of the theorem we describe the connecting
homomorphism $\nu$ in terms of crossed extensions in $\SG$, see
Section \ref{nu}.

It follows from the relationship between $\Sh^3(R,M)$ and crossed
extensions of rings that $\Sh^3(R,M)$ describes homotopy types of
those chain algebras $C_*$ with $H_0(C_*)=R$, $H_1(C_*)=M$, and
$H_i(C_*)=0$ for $i\ne0,1$. On the other hand, it follows from the
relationship between Mac~Lane and topological Hochschild cohomology
that $\HML^3(R,M)$ describes homotopy types of ring spectra
$\Lambda$ with $\pi_i(\Lambda)=0$, $i\ne0,1$, $\pi_0(\Lambda)=R$ and
$\pi_1(\Lambda)=M$ \cite{laza}. Thus our result shows that crossed
extensions of $R$ by $M$ in $\SG$ are algebraic models of such ring
spectra. It follows that the homomorphism $\nu$ is an obstruction
for such a ring spectrum to be representable by a chain algebra.

In Section \ref{obs2cat}
we give an application of our results to
the theory of 2-categories.

\section{Crossed extensions}\label{crex}

We shall apply the following general notion of crossed extension to
the monoidal categories $(\Ab,\t)$, $(\SG,\square)$ and $(\SG,\tl)$
where the category $\SG$ of square groups is defined in \ref{sqgr}.
\subsection{Crossed extensions}\label{crexte}
Let $(\bf V,\boxtimes)$ be a monoidal category and let $L$ be a
monoid in $\bf V$. Recall that a \emph{$L$-biobject} is a tuple
$(A,l,r)$, where $A$ is an object in $\bf V$ and $l:L\boxtimes A\to
A$ and $r:A\boxtimes L\to A$ are respectively left and right actions
of $L$ on $A$ which are compatible in a natural way. We let
$_{L}{\bf V}_{L}$ be the category of $L$-biobjects in $\bf V$. In
particular the monoid structure on $L$ defines also a structure of a
$L$-biobject on $L$. In what follows we always consider $L$ as a
biobject with this particular structure.

A \emph{crossed $L$-biobject} is a diagram $C=(\d:B\to L)$ in the
category $_{L}{\bf V}_{L}$ such that the following diagram commutes:
$$\xymatrix{B\boxtimes B\ar[r]^{\Id\boxtimes \d}\ar[d]_{\d\boxtimes \Id}&B\boxtimes
L\ar[d]^r\\ L\boxtimes B\ar[r]_l&L}$$

Let $R$ be a monoid in $(\V,\boxtimes)$, let $M$ be an $R$-biobject
and assume that exact sequences are defined in $\V$. Then a
\emph{crossed extension} of $R$ by $M$ in $(\V,\boxtimes)$ is an
exact sequence
\begin{equation}\label{crext}
0\to M\xto\iota C_1\xto\d C_0\xto q R\to0
\end{equation}
where $\d$ is a crossed $C_0$-biobject as above, $q$ is a morphism
of monoids and $\iota$ is a morphism in $_{C_0}\V_{C_0}$. A morphism
between crossed extensions of $R$ by $M$ is a commutative diagram
$$
\xymatrix{
0\ar[r]&M\ar[r]\ar[d]_{\Id}&C_1\ar[r]^\d\ar[d]^{f_1}&C_0\ar[r]\ar[d]^{f_0}&R\ar[r]\ar[d]^\Id&0\\
0\ar[r]&M\ar[r]&C_1'\ar[r]^{\d'}&C_0\ar[r]&R\ar[r]&0 }
$$
where $f_0$ is a morphism of monoids and $f_1$ is
$f_0$-biequivariant. Let
\begin{equation}
\xext(R;M)^{\V,\boxtimes}
\end{equation}
be the category of such crossed extensions and morphisms and let
\begin{equation}
\Xext(R;M)^{\V,\boxtimes}
\end{equation}
be the set of connected components of this category.

One readily checks that crossed extensions in $(\Ab,\ox)$ are the
extensions defined in Section \ref{intro}. Hence if $R$ is a ring
and $M$ is an $R$-bimodule then one has canonical bijections (see
\cite{shukla}, \cite[page 42]{jll}, \cite{baumin})
\begin{equation}\label{shuxt}
\Xext(R;M)^{\Ab,\ox}\approx\Sh^3(R;M)
\end{equation}
and
\begin{equation}
\Xext_\Z(R;M)^{\Ab,\ox}\approx\H^3(R;M).
\end{equation}
Here $\Xext_\Z(R;M)^{\Ab,\ox}$ is the set of connected components of
the following subcategory $\xext_\Z(R;M)^{\Ab,\ox}$ of
$\xext(R;M)^{\Ab,\ox}$: its objects are $\Z$-split crossed
extensions $(\iota,\d,q)$ in $(\Ab,\ox)$, that is, with arrows
$\iota$, $\d$ and $q$ admitting a $\Z$-splitting; morphisms in
$\xext_\Z(R;M)^{\Ab,\ox}$ are morphisms $(f_0,f_1)$ in
$\xext(R;M)^{\Ab,\ox}$ such that both $f_0$ and $f_1$ are
$\Z$-split.

\subsection{Square groups}\label{sqgr}
A \emph{square group} is a diagram
$$
A = (\xymatrix{A_e \ar[r]^H & A_{ee} \ar[r]^P & A_e})
$$
where $A_{ee}$ is an abelian group and  $A_e$ is a group. Both
groups are written additively. Moreover $P$ is a homomorphism and
$H$ is a quadratic map, meaning that the \emph{cross effect}
$$(x \mid y)_H = H(x+y)-H(y) -H(x)$$
is linear in $x,y \in A_e$. In addition the following identities are
satisfied
\begin{align*}
(Pa \mid y)_H&=0,\\
P(x\mid y)_H&=[x,y],\\
PHP(a)&=P(a)+P(a).
\end{align*}
Here $[x,y]=-y-x+y+x$, $a,b \in A_{ee}$ and $x,y \in A_e$. It
follows from the first two identities that $P$ maps to the center of
$A_e$. The second identity shows also that $$A^{\mathrm{ad}}:={\sf
Coker}(P)$$ is abelian. Hence $A_e$ is a group of nilpotence class
2. It follows from the axioms that the function $T=HP-\Id_{A_{ee}}$
is an automorphism of $A_{ee}$ and $T^2=\Id_{A_{ee}}$. Moreover, the
function $\Delta:A_e\to A_{ee}$ is linear, where
$$\Delta (x) = HPH(x) + H(x+x) -4H(x)$$
and furthermore  one has the induced homomorphisms $$(-,-)_H:A\ad\t A\ad\to
A_{ee}$$ and $$\Delta:A\ad\to A_{ee}.$$

We refer to \cite{square} and \cite{bjp} for more information on
square groups. We denote by $\SG$ the category of square groups. In
what follows we identify abelian groups and square groups with
$A_{ee}=0$. In this way we obtain a full embedding of categories
$$\Ab\subset\SG
$$
This inclusion corresponds to the fact that any linear functor is
quadratic. The inclusion $\Ab\subset\SG$ has a left adjoint given by
$A\mapsto A\ad$.

The category $\SG$ has two monoidal structures $\square:\SG\x \SG\to
\SG$ \cite{square} and  $\tl:\SG\x \SG\to \SG$ \cite{bjp}, which are
related via a binatural transformation $$\sigma_{X,Y}:X\square Y \to
X\tl Y$$ such that the identity functor together with $\sigma$
defines a lax monoidal functor $\Id:(\SG,\tl)\to (\SG,\square)$
\cite{bjp}. The monoidal category structure $\square$ is  highly
nonsymmetric, while
the monoidal category structure $\tl$ is symmetric.  For the definitions of the products
$\square$ and $\tl$ on $\SG$ we refer the reader to \cite{square} and
\cite{bjp} respectively. Below we shall, however, describe
explicitly the notion of crossed extension in $(\SG,\square)$ and in
$(\SG,\tl)$. Since $(\Ab,\t)$ is a monoidal subcategory both in
$(\SG,\square)$ and in $(\SG,\tl)$, we see that a ring $R$, i.~e. a
monoid in $(\Ab,\t)$, is also a monoid in $(\SG,\square)$ and in
$(\SG,\tl)$. Let $M$ be an $R$-bimodule. Then crossed extensions
$$
0\to M\xto\iota C_{(1)}\xto\d C_{(0)}\xto q R\to0
$$
are defined in $(\SG,\square)$ and in $(\SG,\tl)$ by \eqref{crext}. Such
an extension is \emph{linearly generated} if $R$ as an additive group is
generated by the image of the linear
elements of ${C_{(0)}}_e$ in $R$. Here an element $x\in {C_{(0)}}_e$ is
\emph{linear} provided $H(x)=0$.

As a main result we prove the quadratic analogue of \eqref{shuxt}.

\begin{The}\label{xext}
Let $R$ be a ring and let $M$ be an $R$-bimodule. Then there are
natural bijections
$$
\Xext_L(R;M)^{\SG,\tl}\approx\Xext_L(R;M)^{\SG,\square}\approx\HML^3(R;M)
$$
where the index $L$ indicates the full subcategories of linearly
generated crossed extensions. The first bijection is induced by the
lax monoidal functor $(\SG,\square)\to(\SG,\tl)$.
\end{The}

The proof of this result is given in Section \ref{proof261}.

{\bf Remark}. For the second bijection in the theorem we use the isomorphism
$$
\HML^3(R;M)\cong H^3(\modr R;D_M)
$$
where $\modr R$ is the category of finitely generated  free right $R$-modules
and $D_M=\Hom_R(-,-\ox_RM)$. Here we use the following
interpretation of crossed biobjects from \ref{crexte}.

Let $(\V,\boxtimes)$ be a monoidal category and assume that finite
colimits exist in $\V$. Let $\pair(\V)$ be the category of pairs in
$\V$, objects are morphisms $V=\left(V_1\xto\d V_0\right)$ in $\V$
and morphisms $V\to W$ are pairs $\alpha=(\alpha_1:V_1\to
W_1,\alpha_0:V_0\to W_0)$ in $\V$ with $\d\alpha_1=\alpha_0\d$. Then
$(\pair(\V),\ubx)$ is a monoidal category with $\ubx$ defined by the
diagram with the inner square pushout
$$
\xymatrix{ V_1\boxtimes W_1\ar@{}[ddrr]|{\mathrm{push}}
\ar[rr]^{1 \boxtimes\d}\ar[dd]_{\d \boxtimes 1}
&&V_1\boxtimes W_0\ar[dd]\ar@/^/[dddr]^{\d \boxtimes 1}\\
\\
V_0\boxtimes W_1\ar[rr]\ar@/_/[drrr]_{1 \boxtimes\d}&&(V\ubx W)_1\ar[dr]^\d\\
&&&**[r]V_0 \boxtimes W_0:=(V\ubx W)_0.}
$$
One readily checks that a crossed $L$-biobject $C=(\d:B\to L)$ in
\ref{crexte} is the same as a monoid in $(\pair(\V),\ubx)$. Hence
the action of the monoid $C$ on an object $X$ in $\pair(\V)$ is
defined, compare section 5.1 in \cite{Ba}. In this case we call $X$
a \emph{$C$-module}.

{\bf Addendum}. For a crossed extension $C$ of $R$ by $M$ in either
$(\SG,\square)$ or $(\SG,\tl)$ let $\modr C$ be the category of
finitely generated free left $C$-modules. Then $\modr C$ is a linear
track extension which represents an element
$$
\brk{\modr C}\in H^3(\modr R;D_M)
$$
and the bijections from \ref{xext} carry the component of the
crossed extension $C$ in $(\SG,\tl)$, resp. in $(\SG,\square)$, to the
class $\brk{\modr C}$.

Since a crossed extension $C$ is also a monoid in the category of
pairs, we will also call $C$ a ``pair algebra''. The addendum makes use
of modules over such pair algebras.

\subsection{Square rings and quadratic rings}

A monoid in the
monoidal category $(\SG,\square)$ is termed a \emph{square ring},
while a monoid in the monoidal category $(\SG,\tl)$ is termed a
\emph{quadratic ring}.

More explicitly (see \cite{BHP},
\cite{square}, \cite{iwa}), to provide a square group $Q$ with a
\emph{square ring structure} is the same as to give additionally a
multiplicative monoid structure on $Q_e$. The multiplicative unit of
$Q_e$ is denoted by $1$. One requires that this monoid structure
induces a ring structure on the abelian group $Q\ad$ through the
canonical projection
$$
Q_e\to Q\ad, \ \ \ a\mapsto \bar{a}.
$$
Moreover the abelian group $Q_{ee}$ must be a $Q\ad\t Q\ad\t
(Q\ad)\op$-module with action denoted by $(\bar{x}\t \bar{y})\cdot
a\cdot\bar{z}\in Q_{ee}$ for $\bar{x}, \bar{y}, \bar{z}\in Q\ad$,
$a\in Q_{ee}$. In addition the following conditions must be
satisfied where $H(2)=H(1+1)$
\begin{enumerate}
\item[(i)]
$x(y+z)=xy+xz$

\item[(ii)]
$(x+y)z=xz+yz+P((\bar{x}\t \bar{y})\cdot H(z))$

\item[(iii)]
$(x\mid y )_H=(\bar{y}\t \bar{x})\cdot H(2)$

\item [(iv)]
$T((\bar{x}\t \bar{y})\cdot a\cdot\bar{z})=(\bar{y}\t \bar{x})\cdot
T(a)\cdot\bar{z}$

\item[(v)]
$P(a\cdot x)=P(a)x$

\item[(vi)]
$P((\bar{x}\t \bar{x})\cdot a) =xP(a)$

\item[(vii)]
$H(xy)=(\bar{x}\t \bar{x})\cdot H(\bar{y})+H(x)\cdot \bar{y}$

\end{enumerate}

Under the equivalence ${\sf Quad}(\Gr)\cong\SG$ square rings
correspond to monads on the category of groups, whose underlying
functors lie in ${\sf Quad}(\Gr)$.
A \emph{quadratic ring structure} on a square group $C$ is given by a
multiplicative monoid structure on $C_e$ and a ring structure on
$C_{ee}$. The multiplicative unit of $C_e$ is denoted by $1$. One
requires that these structures satisfy the following additional
conditions.
\begin{enumerate}
\item[(i)] $x(y+z)=xy+xz$,
\item[(ii)] $(x+y)z=xz+yz+P(\cref HyxH(z))$.
\end{enumerate}
Thus $C\ad$ is a ring. Moreover the maps
\begin{align*}
-T&:C_{ee}\to C_{ee},\\
(-\mid -)_H&:C\ad\t C\ad\to C_{ee}
\end{align*}
are ring homomorphisms, in other words one has
\begin{enumerate}
\item[(iii)] $(x\mid y)_H (u\mid v)_H=(xu\mid yv)_H$,
\item[(iv)] $T(ab)+T(a)T(b)=0$.
\end{enumerate}
Let us observe that $T(abc)=T(a)T(b)T(c)$. Furthermore the following equations hold
\begin{enumerate}
\item[(v)] $ P(a\Delta(x))=P(a)x,$
\item[(vi)] $P((x\mid x)_Ha)=xP(a),$
\item[(vii)] $H(xy)=(x\mid x)_HH(y)+H(x)\Delta(y).$
\end{enumerate}

It follows from the axioms that $\Delta:C\ad\to C_{ee}$ is a ring
homomorphism \cite{bjp}.

Let $\QR$ (resp. $\SR$) denote the category of quadratic (resp.
square) rings.  One has the full embedding of categories ${\sf
Rings}\subset\QR $ (resp. ${\sf Rings}\subset\SR $) which identifies
rings with quadratic (resp. square) rings $C$ satisfying $C_{ee}=0$.
This inclusion has a left adjoint given by $R\mapsto R\ad$.

There is also a functor
$$
U:\QR\to \SR
$$
which assigns to a quadratic ring $C$  a square ring, whose
underline square group is the same, while the $C\ad\t C\ad\t
(C\ad)\op$-module structure on $C_{ee}$ is given by
$$ (\bar{x}\t \bar{y})a\bar{z}=(y\mid x)_Ha\Delta(z).$$

The initial object in the category of quadratic rings (resp. square
rings) is $\Znil$, which is given by $$(\Znil)_e=\Z=\ (\Znil)_{ee},
\ \ P=0, \ \ H(x)=\frac{x(x-1)}{2}.$$

We now extend the monoid ring construction to quadratic rings and
square rings. For a monoid $S$ one puts
$$
\Znil[S]_{ee}=\Z[S]\t\Z[S],
$$
where $\Z[S]$ is the free abelian group generated by $S$. We take
$\Znil[S]_e$ to be the free nil$_2$-group generated by $S$. The
homomorphism $P$ is given by $P(s\t t)=[t,s]$, $s,t\in S$, while the
quadratic map $H$ is uniquely defined by
$$
H(s)=0, \ \ (s\mid t)_H=t\t s\ \ s,t\in S.
$$
One has
$$
\Znil[S]\ad=\Z[S].
$$
There is a unique quadratic (resp. square) ring structure on
$\Znil[S]$ for which the multiplication on $\Znil[S]_e$ extends the
multiplication on the monoid $S$ and such that the ring structure
(resp. $\Z[S]\t \Z[S]\t\Z[S]\op$-module structure) on
$\Znil[S]_{ee}=\Z[S]\t\Z[S]$ is the obvious one (resp. is given by
$(x\t y)(s\t t) z=xsz\t ytz$). In this case $Q\ad=\Z[S]$ is the
usual monoid ring of $S$. The functor $\Znil[-]:{\sf Monoids}\to
\QR$ (resp. $\Znil[-]:{\sf Monoids}\to \SR$) is left adjoint to the
functor
$$
{\sf L}:\QR\to \sf Monoids \ \ \ ({\rm resp}. \ \ {\sf L}:\SR\to \sf
Monoids),
$$
where ${\sf L}(Q)$ consists of \emph{linear elements} of $Q$, that
is
$$
{\sf L}(Q)=\{x\in Q_e\mid H(x)=0\}.
$$
The equality $H(xy)=(\bar{x}\t\bar{x})H(y)+H(x)\bar{y}$ shows that
linear elements indeed form a multiplicative submonoid of $Q_e$.

\subsection{Quadratic pair algebras and crossed
square rings}
Now we consider  crossed biobjects in the monoidal
categories $(\SG,\square)$ and $(\SG,\tl)$. Actually we restrict
ourselves to considering only those crossed biobjects $\d:C_1\to
C_0$ which induce isomorphism on ee-level. This is the condition
which implies that $\Cok(\d)$ is a usual  ring and $\Ker(\d)$ is a
usual bimodule. This forces us to introduce the following
definition.

A \emph{quadratic pair module} (qpm for short) is
a morphism of square groups  $\d:C_1\to C_0$
such that the homomorphism $\d_{ee}:{C_1}_{ee}\to{C_0}_{ee}$ is an identity map.
Thus explicitly a quadratic pair module $C$ is
given by a diagram
$$
\xymatrix{&C_{ee}\ar[dl]_P\\
C_{1}\ar[r]_\d&C_{0}\ar[u]_H}
$$
where $C_{1}$ and $C_{0}$ are groups, $C_{ee}$ is an abelian group,
$P$ and $\d$
are group homomorphisms and $H$ is a quadratic map,
and moreover the following
identities are satisfied for any $a\in C_{ee}$, $r,s\in C_{1}$ and
$x,y\in C_{0}$:
\begin{align*}
PH\d P(a)&=2P(a);\\
H(x+\d P(a))&=H(x)+H\d P(a);\\
PH(\d(r)+\d(s))&=PH\d(r)+PH\d(s)+[r,s];\\
\d PH(x+y)&=\d PH(x)+\d PH(y)+[x,y].
\end{align*}
The category of qpm's is denoted by $\QPM$. If $C$ is a qpm, then
$\Im(\d)$ is a normal subgroup of $C_0$ containing the commutator
subgroup of $C_0$. Thus
$$h_0(C):=\Cok(\d)$$
is an abelian group. Moreover
$$h_1(C):=\Ker(\d)$$
is a central subgroup of $C_{1}$.
We have an exact sequence of square groups
$$0\to h_1(C)\to C_{(1)}\xto{(\d,\Id)} C_{(0)}\to h_0(C)\to 0.$$
Here $(C_{(1)})_{ee}=(C_{(0)})_{ee}=C_{ee}$, $(C_{(1)})_{e}=C_{1}$ and
$(C_{(0)})_e=C_{0}.$ The
structural maps are given by
$$ P^{{C}_{(0)}}=\d P, \ \ \ P^{{C}_{(1)}}=P,$$
$$ H^{{C}_{(0)}}=H, \ \ \ H^{{C_{(1)}}}=H\d.$$

A qpm together with crossed biobject structure in the monoidal
category $(\SG,\tl)$ is called a \emph{quadratic pair algebra}
(shortly qpa). Thus a qpa is a  qpm $C$ together with a ring
structure on $C_{ee}$ and a quadratic ring structure on $C_{(0)}$.
Additionally a two-sided action of $C_{0}$ on $C_{1}$ is given, which is
associative and unital and
the following identities are satisfied for all $x,y\in C_{0}$,
$r,s\in C_{1}$, $a,b\in C_{ee}$:
\begin{enumerate}
\item[(i)] $ P((x|x)_Ha\Delta(y))=xP(a)y$
\item[(ii)] $\d(x r y)=x\d (r) y$
\item[(iii)] $  \d(r)s=r\d(s)$
\item[(iv)] $x(r+s)=xr+xs $
\item[(v)] $r(x+y)=rx+ry$
\item[(vi)] $(x+y)r=xr+yr+ P ((y|x)_HH\d(r))$
\item[(vii)] $(r+s)x=rx+sx+P( (s|r)_{H\d}H(x)).$
\end{enumerate}

If $C$ is a qpa, then $C_{(0)}$ is a quadratic ring and the
multiplication on $C_{0}$ yields the multiplication on $h_0(C)$
which equips $h_0(C)$  with a structure of a ring.  Moreover
$h_1(C)$ is a bimodule over $h_0(C)$.

A qpm together with a crossed biobject structure in the monoidal
category $(\SG,\square)$ is called a \emph{crossed square ring}
(shortly csr). Thus a csr is a  qpm $C$ together with a square ring
structure on $C_{(0)}$ and  a two-sided action of $C_{0}$ on $C_{1}$,
which is associative and unital and
such that
 the following identities hold for all $x,y,z,t\in C_{0}$,
$r,s\in C_{1}$, $a,b\in C_{ee}$:
\begin{enumerate}
\item[(i)] $P ((\bar{x}\t \bar{x})\cdot a\cdot y)=x\cdot P(a)\cdot y$
\item[(ii)] $\d(x\cdot r\cdot y)=x\cdot \d (r)\cdot y$
\item[(iii)] $  \d(r)s=r\d(s)$
\item[(iv)] $x(r+s)=xr+xs $
\item[(v)] $r(x+y)=rx+ry$
\item[(vi)] $(x+y) r=xr+yr+ P ((\bar{x}\t \bar{y})\cdot H\d r)$
\item[(vii)] $(r+s)x=rx+sx+P((\bar{\d r}\t \bar{\d s})\cdot Hx).$
\end{enumerate}
In a  crossed square ring  the quotient $R=\Cok(\d)$ has a ring
structure and $\Ker(\d)$ is a bimodule over $R$.
We denote by
$\CS$
the category of  crossed square rings.

Let $R$ be a ring and $M$ be a bimodule over $R$. A \emph{quadratic
ring  extension}  (resp. \emph{crossed square ring extension}) of
$R$ by $M$  is an exact sequence
$$0\to M\xto{i} C_{(1)}\xto{(\d,\Id)} C_{(0)}\xto{p} R\to 0$$
where  $(\d,\Id):C_{(1)}\to C_{(0)}$ is a qpa (resp. csr), the induced
homomorphisms $p:\Cok(\d)\to R$  is an isomorphism of rings and the
induced homomorphism $i:M\to \Ker(\d)$  is an isomorphism of
bimodules over $\Cok(p)$. Here $M$ is considered as a bimodule over
$\Cok(\d)$ via the isomorphism $p:\Cok(\d)\to R$.

Using the definition of the products $\square$ and $\tl$ from
\cite{square}, resp. \cite{bjp}, one readily checks:

\begin{Le}
A crossed extension of $R$ by $M$ in $(\SG,\square)$ is isomorphic
to a crossed square ring extension. A crossed extension of $R$ by
$M$ in $(\SG,\tl)$ is isomorphic to a quadratic ring extension.
\end{Le}

Hence we have explicitly described the objects in the category of
the main theorem \ref{xext}.
\begin{Exm}\label{ztildanil} Let $Q$ be a square ring. One can
consider the quotient $Q_{ee}/(\Id-T)$, where as usual $T=HP-\Id$.
Let $\tP:Q_{ee}\to Q_{ee}/(\Id-T) $ be the canonical projection. It
is clear that the homomorphism $P:Q_{ee}\to Q_e$ factors through
$Q_{ee}/(\Id-T)$. We denote by $\d:Q_{ee}/(\Id-T)\to Q_e$ the
quotient map. Then $$\xymatrix{& Q_{ee}\ar[dl]_{\tP}\\
Q_{ee}/(\Id-T)\ar[r]_{\ \ \d} &Q_e\ar[u]_H}$$ is a crossed square
ring. Thus Theorem \ref{xext} assigns to any square ring $Q$ an
element in $\HML^3(Q\ad, Q^{re})$, where $Q^{re}
=\Ker(\d:Q/(\Id-T)\to Q_e)$. In particular, for the square ring
$\Znil$ one obtains the following crossed square ring
$$\xymatrix{& \Z\ar[dl]_{1}\\
\Z/2\Z\ar[r]_{0}&\Z\ar[u]_H}$$ where $H(x)=\frac{x^2-x}{2}$ which
defines an element of $\HML^3(\Z,\Z/2\Z)=\Z/2\Z$, which is actually
the  generator.

\end{Exm}

\subsection{The homomorphism $\nu:\HML^3(R,M)\to
\H^0(R,\,_2M)$}\label{nu}
Let $R$ be a ring and let $M$ be a bimodule over $R$. Take a crossed
square ring extension $(\d)$ of $R$ by $M$
$$\xymatrix{
0\ar[r]& M\ar[r]^i& C_{(1)}\ar[r] ^{(\d,\Id)} &C_{(0)} \ar[r]^p&R\ar[r]&0}$$ and set
$$\upsilon(w):=P H(2)$$
Since $H(1)=0$ it follows that $\d P H(2)=\d P(1\mid 1)_H=0$. On the
other hand $2P H(2)=P H\d P H(2)=0$. Thus $\upsilon(w)\in \, _2\,M$.
Actually $$\upsilon(w)\in \, \H^0(R,\, _2\,M)$$ and $\nu$ yields a
well-defined map

$$
\nu:\Xext_L(R;M)^{\SG,\square}\to\H^0(R;{}_2M).
$$

\begin{Le}
Kernel of $\nu$ coincides with the image of
$$
\Xext(R;M)^{\Ab,\t}\to\Xext_L(R;M)^{\SG,\square}.
$$
\end{Le}

In fact we obtain the lemma directly by the exact sequence in
section 1 and the bijections (\ref{shuxt}) and \ref{xext}. We snow show
the lemma more directly in terms of crossed extensions.

\begin{proof}
If $\d$ is a crossed ring extension then $C_{ee}=0$ and a fortiori
$H=0$, thus $\nu(\d)=0$. Conversely, assume $(\d)$ is a crossed
square ring extension with $\nu(\d)=0$. Without loss of generality
one can assume that $C_{(0)}$ is a monoid square ring
$C_{(0)}=\Znil[L]$ (see Section \ref{proof261} below). In this case
$C_{ee}=\Z[L]\t \Z[L]$ and $H(2)=(1\mid 1)_H=1\t 1.$ The equality
$P((\bar{x}\t \bar{x}\cdot m\cdot y)=x\cdot P(m)\cdot y$ shows that
$P$ factors through $\Lambda^2(\Z[L])$. Thus one gets the following
diagram
$$\xymatrix{&&& 0\ar[d]&&\\
&& \Lambda ^2(\Z[L])\ar[r]^\id\ar[d]^{\tP}&\Lambda ^2(\Z[L])\ar[d]^{[-,-]}&&\\
0\ar[r]&M\ar[r]\ar[d]_\Id&C_1\ar[r]^\d\ar[d]&C_{0}\ar[r]\ar[d]&R\ar[r]
\ar[d]_\Id&0\\
0\ar[r]&M\ar[r]&\Cok{\tP}\ar[r]\ar[d]&\Z[L]\ar[r]\ar[d]&R\ar[r]&0\\
&&0&0&& }$$ Since $C_{0}$ is a free nil$_2$-group on $L$ the
commutator map $[-,-]$ is a monomorphism. It follows that $\tP$ is
also a monomorphism, the bottom row is exact  and $(\d)$ is
equivalent to
$$0\to M\to \Cok{\tP}\to \Z[L]\to R\to 0$$
which is a crossed ring extension in $(\Ab,\t)$.
\end{proof}
\subsection{Application to ring spectra} Since Mac~Lane
cohomology and topological Hochschild cohomology are isomorphic for
discrete rings it follows that for any ring spectrum $\Lambda$ with
$\pi_i(\Lambda)=0$ for $i\ne0,1$ there is a well-defined element
$k(\Lambda)\in \HML^3(\pi_0(\Lambda),\pi_1(\Lambda))$ known as the
first Postnikov invariant (see \cite{laza}) and any element in this
group comes in this way. Thus linearly generated crossed square
rings and quadratic pair algebras can be used to model such ring
spectra. The explicit functor from the category of crossed square
rings to the category of ring spectra can be constructed as follows.
By Corollary \ref{lodjan} below one can associate to any crossed
square ring an internal groupoid in the category of square rings and
hence an internal groupoid in the category of algebraic theories
(see \ref{sqringtheories} below). Now using the nerve construction one
obtains a simplicial object in the category of algebraic theories.
Then one can use the well-known construction of Schwede
\cite{schwede} to obtain a ring spectrum in a functorial way.

\section{Recollections}
\subsection{Preliminaries on double categories and internal categories}

Let $\A$ be a category with finite limits. An \emph{internal
category} $C$ in $\A$ consists of the following data: objects $C_0$
(\emph{object of objects}), $C_1$ (\emph{object of morphisms}) and
morphisms  $s,t:C_1\to C_0$ (\emph{source} and \emph{target}),
$i:\C_0\to \C_1$ (\emph{identity}), $m:C_2\to C_1$
(\emph{composition}) satisfying associativity and unitality
conditions. Here $C_2$ is defined by the pullback diagram
$$
\xymatrix{
G_2\ar[r]^{p_1}\ar[d]_{p_2}& G_1\ar[d]^s\\
G_1\ar[r]_t&G_0
}
$$
We denote by $\Cat(\A)$ the category of internal categories in $\A$.
Let us also recall that an internal category $C$ is called an
\emph{internal groupoid} provided the diagram
$$
\xymatrix{
G_2\ar[r]^{m}\ar[d]_{p_2}& G_1\ar[d]^s\\
G_1\ar[r]_s&G_0
}
$$
is a pullback diagram. We denote by $\Gpd(\A)$ the category of
internal groupoids in $\A$.

An internal category in the category of sets $\Sets$ is nothing but
a small category, while a groupoid object in the category of sets
$\Sets$ is a groupoid. We write $\Cat$ and $\Gpd$ instead of
$\Cat(\Sets)$ and $\Gpd(\Sets)$.

Let $A$ be an object of $\A$, then we can consider the internal
groupoid $A^{dis}$ with $(A^{dis})_0=A=(A^{dis})_1$ and $s=t=\Id_A$.
An internal category is called \emph{discrete} if it is isomorphic
to $A^{dis}$ for some $A$. We will need also an internal groupoid
$A^{adis}$ with $(A^{adis})_0=A$, $(A^{adis})_1=A\x A$, where $s$
and $t$ are the projections. An internal category is called
\emph{antidiscrete} if it is isomorphic to $A^{adis}$ for some $A$.

Let $\B$ be a category with finite limits and let $F:\A\to\B$ be a
functor which preserves finite limits. Then obviously $F$ yields
functors $\Cat(\A)\to\Cat(\B)$ and $\Gpd(\A)\to\Gpd(\B)$ which will
be also denoted by $F$.

Let us recall that a \emph{double category} is an internal category
in the category $\Cat$ of small categories. Let $\D$ be a double
category with the object category $\D_0$ and morphism category
$\D_1$.

We have a functor ${\sf Ob}:\Cat\to\Sets$, which assigns to a
category $\C$ the set of objects of $\C$. Since ${\sf Ob}$ preserves
inverse limits, for any double category $\D$ we obtain a category
$O(\D)$, whose morphisms are objects of $\D_1$ and objects are
objects of $\D_0$. A double category $\D$ is a \emph{2-category} if
$O(\D)$ is a discrete category. Equivalently a 2-category is a
category enriched in the category $\Cat$. Let us recall how one gets
such an enrichment.

Let $\D$ is a 2-category. Then objects of the category $\D_0$ are
called simply objects of $\D$, while morphisms of the category
$\D_0$ are called simply morphisms of $\D$. Let $f,g:A\to B$ be
morphisms of $\D$. Then $A$ and $B$ are also objects in $\D_1$ and
we can consider the set of all morphisms $\alpha:A\to B$ in $\D_1$
such that $s(\alpha)=f$ and $t(\alpha)=g$. Such an $\alpha$ is
called a 2-morphism from $f$ to $g$. Thus for objects $A$ and $B$ we
have a category $\D(A,B)$ with objects morphisms from $A$ to $B$ in
the category $\D_0$ and morphisms from $f:A\to B$ to $g:A\to B$
being all 2-morphisms from $f$ to $g$.

Conversely, if $\B$ is a a
category enriched in the category $\Cat$, then one can consider the
following categories $\B_0$ and $\B_1$. The category $\B_0$ has the
same objects as $\B$, while morphisms in $\B_0$ are 1-arrows of
$\B$. The category $\B_1$ has the same objects as $\B_0$. The
morphisms $A\to B$ in $\B_1$ are 2-arrows $\alpha:f\then f_1$ where
$f,f_1:A\to B$ are 1-arrows in $\B$. Composition in $\B_1$ is given
by $(\bb:x\then x_1)(\aa:f\then f_1):=(\bb\aa:xf\then x_1f_1)$,
where
$$
\bb\aa=\bb f_1+x\aa=x_1\aa+\bb f.
$$
One furthermore has the source and target functors
$$
\xymatrix{
\B_1\ar@<1ex>[r]^s\ar@<-1ex>[r]_t&\B_0
},
$$
with $s(\alpha:f\then f_1)=f$, $t(\alpha:f\then f_1)=f_1$, and the
``identity'' functor $i:\ta_0\to\ta_1$ assigning to an 1-arrow $f$
the identity 2-arrow $0_f:f\then f$. One easily sees that in this
way we obtain a double category such that after applying the functor
${\sf Ob}:\Cat\to\Sets$ one gets a discrete category.

\subsection{Preliminaries on Baues-Wirsching cohomology of small categories
}\label{relcoh}
For a small category $\c$ we denote by $\f\c$ the \emph{category of
factorizations} of $\c$ \cite{BW}. Objects of $\f\c$ are morphisms
of $\c$, and a morphism from $\alpha:x\to y$ to $\beta:u\to v$ is a
pair $(\nu:u\to x,\psi:y\to v)$ of morphisms in $\c$ such that
$\beta = \psi \alpha \nu$, that is, one has a commutative diagram
$$
\xymatrix{ x\ar[rr]^{\alpha}  &  & y  \ar[d]^{\psi}\\
u \ar[u]^{\nu}\ar[rr]^{\beta}  &       &    v      }
$$
Composition in $\f\c$ is defined by $(\nu,\psi)(\nu',\psi')= (\nu'
\nu, \psi \psi')$. A \emph{natural system} on $\c$ is a covariant
functor $D:\f\c\to\Ab$. Now, following \cite{BW}, one defines the
cohomology $\H^*(\c,D)$ as the cohomology of the cochain complex
$\F^*(\c,D)$ given by
$$
\F^n(\c,D)= \prod_{ c_0  \stackrel{\alpha_1}{\leftarrow}  \cdots
 \stackrel{\alpha_n}{\leftarrow}  c_n} D_{ \alpha_1 \cdots \alpha_n}
$$
with the coboundary map
$$
d : \F^n(\c,D) \to \F^{n+1}(\c,\D)
$$
given by
$$
\begin{array}{ll}
(df)(\alpha_1, \cdots, \alpha_{n+1}) &= (\alpha_1)_* f(\alpha_2,
\cdots, \alpha_{n+1})  \\ \\
& + \sum_{i=1}^n (-1)^i f( \alpha_1, \cdots, \alpha_i
\alpha_{i+1}, \cdots, \alpha_{n+1}) \\ \\
& + (-1)^{n+1} (\alpha_{n+1})^* f (\alpha_1, \cdots, \alpha_{n}).
\end{array}
$$
Here, and in the rest of the paper, we use the following convention.
For a diagram $u\xto{\beta}x\xto{\alpha}y\xto{\gamma}v$ and elements
$a \in D_\beta$, $ b \in D_{\gamma}$, we write $\alpha_*a$ and
$\alpha^*b$ for the image of the elements $a$ and $b$ under the
homomorphisms $D(id_u, \alpha): D_\beta \to D_{\alpha\beta}$ and
$D(\alpha, id_v): D_\gamma \to D_{\gamma\alpha}$ respectively.

We also need the relative cohomologies of small categories. Let
$p:\ka\to\c$ be a functor which is identity on objects and
surjective on morphisms. Let $D:\f\c\to\Ab$ be a natural system on
$\c$. We have an induced natural system $p^*D$ on  $\ka$ given by
$g\mapsto D_{pg}$, which we will, abusing notation, still denote by
$D$. Then $p$ yields a monomorphism of cochain complexes
$\F^*({\c},D)\to \F^*(\ka,D)$.  We let $\F^*(\c,\ka;D)$ be the
cokernel of this homomorphism. The {\it $n$-th dimensional relative
cohomology} $\H^n(\c,\ka;D)$ is defined as the $(n-1)$-th homology
of the cochain complex  $\F^*(\c,\ka;D)$. Then one has an exact
sequence
$$0\to \H^0({\c},D)\to \H^0(\ka,D)\to \H^1(\c, \ka;D)\to \cdots \to$$
$$\to \H^n({\c },D)\to \H^n(\ka,D)\to \H^{n+1}(\c,\ka;D)\to \cdots.$$

We have a functor $\f\c \to \c\op \times \c$ which sends an arrow
$\alpha:c \to d$ to the pair $(c,d)$.  This functor allows us to
conclude that any bifunctor gives rise to a natural system. Thus for
any bifunctor $D:\c\op \times \c\to \Ab$ we have well-defined
cohomology groups $H^*(\c,D)$.
Among many equivalent definitions of the Mac~Lane cohomology
\cite{JP} for our purposes  the most convenient  is via the
Baues-Wirsching cohomology of small categories \cite{BW}. Let $R$ be
a ring. Let $\Mod R$ be the category of  right $R$-modules and let
$\modr R$ be the full subcategory of finitely generated free right
$R$-modules. To avoid set-theoretic complications we will assume
that objects of $\modr R$ are natural numbers and morphisms from $y$
to $x$, $x,y\in \N$  are $(x\x y)$-matrices with entries in $R$. We
write $f=(f_i^k)$ for a morphism $y\to x$, where $f_i^k\in R, 1\le
i\le x, 1\le k\le y$.

For an $R$-$R$-bimodule $M$, we denote by $D_M:(\Modr R)\op\x \Modr
R\to \Ab$ the bifunctor given by
$$
D_M(X,Y):=\Hom_R(X,Y\t_R M), \ \
X, Y \in \Modr R.
$$
Now one defines the \emph{Mac~Lane cohomology of $R$ with
coefficients in $M$} by
$$
\HML^*(R,M):=\H^*(\modr R,D_M).
$$
We refer to \cite{JP} and Chapter 13 of \cite{jll} for relationship
between different definitions of Mac~Lane cohomology. We use this
definition of $\HML^*$ also if $R$ is a square ring or a quadratic
ring.

\subsection{Third Baues-Wirsching cohomology and linear track extensions}

We recall the relationship between third Baues-Wirsching cohomology
and linear track extensions. We start with recalling the definition of track
categories.

A \emph{track category} is a groupoid enriched category, i.~e. a
2-category such that all of its 2-morphisms are invertible.
Equivalently a track category $\ta$ is an internal groupoid in the
category $\Cat$ such that $O(\ta)$ is a discrete category. We will
use the following notation for track categories. Composition of
morphisms will be denoted by juxtaposition; for 2-arrows we will use
additive notation, so composition is + and identity 2-arrows are
denoted by 0. The hom-category for objects $A$, $B$ of a track
category will be denoted by $\hog{A,B}$. If there is a 2-arrow
$\alpha:f\then g$ between maps $f,g\in\Ob(\hog{A,B})$, we will say
that $f$ and $g$ are homotopic and write $f\ho g$. We have the
\emph{homotopy category} $\ta_\ho=\ta_0/\ho$. Objects of $\ta_\ho$
are  objects in $\Ob(\ta)$, while morphisms of $\ta_\ho$ are
homotopy classes of morphisms in $\ta_0$. A map $f$ in $\ta$ is
called a \emph{homotopy equivalence} if the class of $f$ in
$\ta_\ho$ is an isomorphism.

Two track categories $\ta$, $\ta'$ are called \emph{weakly
equivalent} if there is an enriched functor $F:\ta\to\ta'$ which
induces equivalences of hom-groupoids
$\hog{X,Y}_\ta\to\hog{FX,FY}_{\ta'}$ and is \emph{essentially
surjective}, i.~e. any object of $\ta'$ is homotopy equivalent to
one of the form $FX$.

 Let $\c$ be a small category and let
$D$ be a natural system on $\c$. A \emph{linear track extension of $\c$ by $D$} denoted by
$$
0\to D\to \ta_1 \rightrightarrows \ta_0\to \c \to 0
$$
is a pair $(\ta, \tau)$. Here $\ta$ is a track category equipped
with a functor $q:\ta_0\to\c$ which is full and identity on objects.
In addition for maps $f,g$ in $\ta_0$ we have $q(f)=q(g)$ iff $f\ho
g$. In other words the functor $q$ identifies $\c$ with $\ta_\ho$.
Furthermore, for each map $f:A\to B$ in $\ta_0$ there is given an
isomorphism of groups $\tau_f:D_{qf}\to {\ta}(f,f)$, such that for
any $\xi:f\then g$ and $a\in D_{qf}=D_{qg}$ one has
$$
\xi+\tau_f(a)=\tau_g(a)+\xi.
$$
Furthermore for any diagram $\buildrel e \over \longrightarrow
\buildrel f \over \longrightarrow \buildrel h \over \longrightarrow$
additionally one has
$$h_*\tau_f(a)=\tau_{hf}(h_*a),$$
$$ e^*\tau_f(a)=\tau_{fe}(e^*a).$$

For a category $\c$ and a natural system $D:\f\c\to \Ab$ we denote
by ${\Track}(\c,D)$ the category of all linear track extensions of
$\c$ by $D$, where the morphisms are the obvious ones.

Linear track extensions of categories were first described in the
preprint of \cite{B} and the following theorem in a slightly
different terminology first was
proved in \cite{P1} (see also \cite{P2}) and was proved
by different methods in \cite{BD}.

\begin{The}\label{h3}\cite{P1} For a small category $\c$ and a natural system
$D:\f\c\to \Ab$ there exists a natural bijection between the set of
connected components of the category $\Track(\c,D)$ and third
cohomology:
$$
\pi_0(\Track(\c,D))\cong \H^3(\c,D).
$$
\end{The}

The proof of Theorem \ref{h3} given in \cite{P1} and \cite{P2} is
based on the following Theorem \ref{fardobiti}, which is going to be
crucial in this paper as well.

Let $p:\ka\to \c$ be a functor which is identity on objects and
surjective on morphism. Let $D:\f\c\to \Ab$ be a natural system on
$\c$. We denote by $\Track(\c,\ka;D)$ the subcategory of
$\Track(\c,D)$ whose objects are track categories $\ta$ satisfying
$\ta_0=\ka$,
$$
0\to D\to \ta_1 \rightrightarrows \ka \buildrel {q}\over \to \c \to
0,
$$
whereas morphisms are those morphisms in $\Track(\c,D)$ which are
identity on $\ka$.

\begin{The}\label{fardobiti}\cite{P1},\cite{P2}
For a small category $\c$, a bifunctor $D:\c\op\x \c\to \Ab$ and a
functor $p:\ka\to \c$ which is identity on objects and surjective on
morphisms, the category $\Track(\c,\ka;D)$ is a groupoid and there
exists a natural bijection
$$
\pi_0(\Track(\c,\ka;D))\cong \H^3(\c,\ka;D).
$$
\end{The}

\subsection{Relative track extensions of algebraic theories}
An \emph{algebraic  theory} is a category with finite coproducts. A
\emph{morphism of algebraic theories} is a functor preserving finite
coproducts. We denote the coproduct by $\vee$. Let $\c$ be an
algebraic theory. A natural system $D:\f\c\to \Ab$ is called
\emph{cartesian} if for any arrow $f:c=c_1\vee \cdots \vee c_n\to d$
the natural map
$$
D_f\to D_{f_1}\x \cdots \x D_{f_n},
$$
given by $x\mapsto ((i_1)^*x,\cdots ,(i_n)^*x)$, is an isomorphism.
Here $i_k:c_k\to c$ is the standard inclusion and $f_k=i_k\circ f:c_k\to
d$. For example if $D:\c\x\c\op\to\Ab$ is a bifunctor such that for
all $c,d$ and $x$ from $\c$ one has an isomorphism
$$
D(c\vee d,x)\cong D(c,x)\x D(d,x)
$$
natural in $c,d$ and $x$, then the natural system corresponding to $D$ is
cartesian.

Let $\c$ be an algebraic theory and let $D:\f\c\to\Ab$ be a
cartesian natural system. A  track extension
$$
0\to D\to \ta_1 \rightrightarrows \ta_0\xto{p} \c \to 0
$$
is called a \emph{track extension of algebraic theories} if $\ta_0$
is an algebraic theory and the functor $p$ is a morphism of
algebraic theories.

\begin{Le}\label{trth} Let $\c$ be an algebraic theory and let $D:\f\c\to\Ab$ be a cartesian
natural system. Let $0\to D\to \ta_1 \rightrightarrows \ta_0\xto{p} \c \to 0$ be a
track extension of algebraic theories. Then $\ta_1$ is also an algebraic theory
and $s,t:\ta_1\to \ta_0$ are morphisms of algebraic theories.
\end{Le}

\begin{proof}
Let $\aa:f\then g$ and $\aa':f'\then g'$ be tracks, where $f,g:A\to
B$ and $f',g': A'\to B$ are 1-morphisms. We have to show that there
is a unique track $(\aa, \aa'):(f,g)\then (f', g')$ such that
$i_A^*(\aa,\aa')= \aa$ and $ i_{A'}^*(\aa, \aa')= \aa'$, where
$(f,g):A\vee A'\to B$ is the unique 1-morphism with $i_A(f,g)=f$ and
$i_{A'}(f,g)=g$. Here $i_A:A\to A\vee B$ is the canonical inclusion
and similarly for $i_{A'}$. First we show the existence of such a
track. By assumption $f \ho g$ and $f'\ho g'$. Since $p$ preserves
finite coproducts, it follows that $(f, f')\ho (g, g')$. Hence there
exists a track $\eta: (f, f')\then (g, g')$. Since
$i_A^*(\eta):f\then g$, there exists a unique element $x\in D_{pf}$
such that $i_A^*(\eta)=\aa+ \sigma_{f}(x)$. Similarly there exists a
unique element $x'\in  D_{pf'}$ such that $i_{A'}^*(\eta)=\aa'+
\sigma_{f'}(x')$. By our assumptions there is a unique element $y\in
D_{(pf,pf'}$ such that $i_A(y)=x$ and $i_{A'}(y)=x'$. Then the track
$\xi=\eta-\sigma_{(f, f')}(y)$ satisfies the condition required. To
prove uniqueness one observes that if $\xi$ and $\eta$ both satisfy
the condition, then they will differ by an element $z\in
D_{(pf,pf')}$, whose restrictions to $D_{pf}$ and $D_{pf'}$ are
zero, hence it is itself zero and the lemma follows.
\end{proof}

\begin{Le}\label{trcogr}
Let $0\to D\to \ta_1 \rightrightarrows \ta_0\xto{p} \c \to 0$ be a
track extension of algebraic theories and let $\nu:X\to X\vee X$ be
an internal cogroup in $\ta_0$. Then $X$ is also a cogroup in
$\ta_1$, where the cogroup structure is given by the morphism
$0:\nu\then\nu$.
\end{Le}

\begin{proof}
By Lemma \ref{trth} the ``identity functor'' $\ta_0\to\ta_1$
respects finite coproducts and therefore carries a cogroups to
cogroups.
\end{proof}

\subsection{Quadratic functors, quadratic categories and square objects}
We now recall the relationship between square groups and quadratic
functors. We consider endofunctors $F:\Gr\to\Gr$ of the category of
groups with $F(0)=0$. Additionally we assume that $F$ preserves
filtered colimits and reflexive coequalizers. The last condition
means that for any simplicial group $G_*$ the canonical homomorphism
$\pi_0(F(G_*))\to F(\pi_0(G_*))$ is an isomorphism. Such a functor
$F$ is completely determined by the restriction of $F$ to the
subcategory of finitely generated free groups.

The second cross-effect  $F(X|Y)$ of $F$ is a bifunctor defined via
the short exact sequence
$$
0\to F(X|Y)\to F(X\vee Y)\to F(X)\x F(Y)\to 0.
$$
Here  $\vee$ denotes the coproduct in the category of groups and the
last map is induced by the canonical projections:
$r_1=(\Id_X,0):X\vee Y\to X$ and $r_2=(0,\Id_Y):X\vee Y\to Y$. A
functor $F$ is called \emph{linear} if the second cross effect
vanishes. Moreover, $F$ is called \emph{quadratic} if $F(X|Y)$ is
linear in $X$ and $Y$. Let ${\sf lin}(\Gr)$ (resp. ${\sf
Quad}(\Gr)$) be the category of such linear (resp. quadratic)
endofunctors. Any endofunctor in ${\sf lin}(\Gr)$ is isomorphic to a
functor $T$ of the form $ T(X)=A\t X\ab$ where $A$ is an abelian
group. Therefore there is an equivalence of categories
$$
 {\sf lin}(\Gr)\simeq \Ab.
$$
Let $F:\Gr\to \Gr$ be a quadratic functor. We associate with $F$ a
square group ${\sf cro}(F)$ as follows. We put
$$
{\sf cro}(F)_e=F(\Z), \ \ {\sf cro}(F)_{ee} =F(\Z\mid\Z).
$$
The homomorphism $P$ of the square group ${\sf cro}(F)$ is the
restriction of the homomorphism $(\Id,\Id)_*: F(\Z \vee \Z)\to
F(\Z)$. We denote by $e_1$ and $e_2$ the canonical free generators
of $\Z\vee\Z$. The map $H$ is given by
$$
H(x)=\mu_*(x)-p_2(\mu_*x)-p_1(\mu_*x)
$$
Here $\mu:\Z\to \Z \vee\Z$ is the unique homomorphism which sends $1$ to
$e_1+e_2$, while $p_1$ and $p_2$ are endomorphisms of $\Z\vee\Z\to \Z\vee\Z$
such that $p_i(e_i)=e_i$, $i=1,2$ and $p_i(e_j)=0$, if $i\not =j$.

The main result of \cite{square} claims that the functor
$$
{\sf cro}: {\sf Quad}(\Gr)\to \SG
$$
is an equivalence of categories. Under this equivalence  square rings
corresponds to monads on the category of groups, whose underlying
functors lie in ${\sf Quad}(\Gr)$.

Let $\C$ be an algebraic theory with zero object $0$. We will say
that $\C$ is equipped with a structure of \emph{quadratic theory} if
each object $C$ in $\C$ is equipped with a cogroup structure
$\nu_C:C\to C\vee C$ and the functor $\C(C,-):\C\to \Gr$ is
quadratic. Thus for all $X$ and $Y$ in $\C$ one has the following
short exact sequence of groups
$$
0\to \C(C;X\mid Y)\to \C(C,X\vee Y)\to \C(C,X)\x \C(C,Y)\to 0
$$
and $\C(C;X\mid Y)$ is linear in $X$ and $Y$. This definition is
equivalent but not identical to the one given in \cite{BHP}.

\begin{Le}\label{lec} Let $\c$ be an additive category,
$D:\c\op\x \c\to \Ab$ be a biadditive bifunctor and
$$0\to D\to \ta_1 \rightrightarrows \ta_0\xto{p} \c \to 0$$
be a linear track extension. If $\ta_0$ is a quadratic theory and
$p$ preserves finite coproducts, then $\ta_1$ is also a quadratic
theory.
\end{Le}

\begin{proof} By Lemma \ref{trth} the category $\ta_1$ is an algebraic theory.
It is quite easy to show that the zero object in $\ta_0$ remains
also a zero object in $\ta_1$. By Lemma \ref{trcogr} any object in
$\ta_1$ has a canonical cogroup structure. We claim that
$$\ta_1(X, Y\mid Z)\cong \ta_0(X, Y\mid Z)\x \ta_0(X, Y\mid Z),$$
which implies that $\ta_1(X, -)$ is a quadratic functor and hence
Lemma. To prove the claim we put
$$
r_1=(\Id,0):Y\vee Z\to Y \  \ \mathrm{ and} \ \ r_2=(0,\id):Y\vee Z\to Z.
$$
By definition of the cross-effect $\ta_0(X, Y\mid Z)$ consists of
1-morphisms $f:C\to Y\vee Z$ such that $r_1f=0=r_2f$. On the other
hand $\ta_1(X, Y\mid Z)$ consists of tracks  $\aa:f\then g$ such
that $r_1f=0=r_2f$, $r_1g=0=r_2g$ and $r_{1*}\aa=0=r_{2*}\aa$. Thus
our claim is equivalent to the following: suppose $f,g:X\to Y\vee Z$
are 1-morphisms such that $r_1f=0=r_2f$, $r_1g=0=r_2g$. Then there
exists a unique track $\aa:f\then g$ such that
$r_{1*}\aa=0=r_{2*}\aa$. If $\aa$ and $\bb$ both satisfy the
assertion, then $\bb=\aa =\sigma_f(x)$ with $x\in D(X,Y\vee
Z)=D(X,Y)\oplus D(X,Z)$. The conditions $r_{1*}\aa=0=r_{2*}\aa=
r_{1*}\bb=r_{2*}\bb$ show that $x=0$. Hence we proved uniqueness.
Now we prove the existence. Since $\c$ is an additive category, $p$
respects coproducts and $r_1f=0=r_2f$ it follows that $p(f)=0$ in
$\c$. Similarly $p(g)=0$. In particular $f\ho g$. Thus there exist a
track $\xi:f\then g$. Then $r_{1*}(\xi):0\then 0$. Hence there
exists a unique $y\in D(X,Y)$ such that $r_{1*}(\xi)=\sigma_0(y)$.
Similarly $r_{2*}(\xi)=\sigma_0(z)$ for uniquely defined $z\in
D(X,Z)$. Since $D$ is biadditive, we have $(y,z)\in D(X,Y\vee Z)$.
It is clear that $\aa=\xi-\sigma_f(y,z)$ satisfies the conditions of
the claim.
\end{proof}

\subsection{Square rings and single sorted quadratic
theories}\label{sqringtheories} We recall the relationship between
square rings and quadratic categories \cite{BHP}. A quadratic theory
is a \emph{single sorted quadratic theory} if the objects of $\C$
are natural numbers and the coproduct on objects corresponds to the
addition of natural numbers. Thus each object ${\bf n}$ in $\C$ is
an $n$-fold coproduct of $\bf 1$. We additionally require that the
cogroup structure on ${\bf n}$ is the $n$-fold coproduct of the
cogroup structure on $\bf 1$.

Assume $\C$ is a single sorted quadratic theory. Then one has the
square ring ${\sf cro}(\C)$ with
$$
{\sf cro}(\C)_e={\sf cro}(\C({\bf 1},-))_e=\C({\bf 1},{\bf 1})
$$
and
$${\sf cro}(\C)_{ee}={\sf cro}(\C({\bf 1},-))_{ee}=\C({\bf 1};{\bf 1}\mid {\bf 1}).$$
The  main result of \cite{BHP} shows that the functor ${\sf cro}$
from the category of single sorted quadratic theories to the
category of square rings is an equivalence of categories. The
inverse functor is given by $Q\mapsto\modr Q$. Here the objects of
the category $\modr Q$ are  natural numbers, while morphisms from
$y$ to $x$, $x,y\in \N$ are defined by product sets
$$
{\sf Mor}(y,x):=(\prod_{k=1}^y\prod_{i=1}^x
Q_e)\x(\prod_{k=1}^y\prod_{1\le i<j\le x}Q_{ee}).
$$
For a morphism $f:y\to x$ we write $f=(f_i^k,f^k_{ij})$.
If $g=(g^s_k,g^s_{kl})$ is a morphism $z\to y$, then the composite
$fg=((fg)_i^s,(fg)^s_{ij})$ is given by
$$(fg)^s_i=f^1_i\circ g_1^s+\cdots +f_i^y\circ
g^s_y+\sum_{k<l}P((\bar{f}_i^k\t \bar{f}_i^l)\cdot g^s_{kl})$$
$$(fg)^s_{ij}=\sum_k(f^k_{ij}\cdot \bar{g}^s_k+
\sum_{i<l}((\bar{f}_i^k\t \bar{f}_i^l)\cdot g^s_{kl}+ (\bar{f}_i^l
\t \bar{f}_j^{k} )\cdot T g^s_{kl}+\overline{(f^l_i\cdot
g^s_l)}\t\overline{(f^k_j\cdot g^s_k)}\cdot H(2))$$ Actually $\modr
Q$ is a single sorted quadratic theory, the group structure on hom's
is defined by the formula:
$$(f_i^k,f_{ij}^k)+({f'}_i^k,{f'}_{ij}^k)=(f_i^k+{f'}_i^k,f_{ij}^k+{f'}^k_{ij}+
e_{ij}^k)$$ where $$e_{ij}^k=(\bar{f}_i^k\t \bar{f'}_j^k)\cdot
H(2).$$ To  get more hints on the category $\modr Q$, we recall that
a \emph{right $Q$-module} \cite{BHP} is nothing but a right
$Q$-object in the monoidal category $(\SG,\square)$.  More
explicitly, a right $Q$-module is a group $M$ together with maps
$M\x Q_{e}\to M$, $(m,x)\to mx$ and $M\x M\x Q_{ee}\to M$,
$(m,n,a)\mapsto [m,n]_a$ satisfying the following identities.
\begin{align*}
m1&=m,\\
(mx)y&=m(xy),\\
m(x+y)&=mx+my,\\
(m+n)x&=mx+nx+[m,n]_{H(a)},\\
mP(a)&=[m,m]_a,\\
[m,n]_{Ta}&=[n,m]_a,\\
[mx,ny]_a&=[m,n]_{(x\otimes y)a},\\
[[m,n]_a,z]_b&=0.
\end{align*}
Moreover $[m,n]_a$ is linear in $m,n$ and $a$ and lies in the center
of $M$. We denote by $\Mod Q$ the category of all right $Q$-modules.
It is a standard fact of universal algebra that the forgetful
functor $\Mod Q\to \sf Sets$ has the left adjoint, whose values on a
set $X$ is called the free right $Q$-module generated by the set
$X$. Now one checks directly \cite{BHP} that the category $\modr Q$
is equivalent to the category of finitely generated free right
$Q$-modules.

Let us observe that for $Q=\Znil$, the category of right
$\Znil$-modules is nothing but the category $\Nil$ of groups of
nilpotence class two. More generally, if $S$ is a monoid then the
category of right modules over the square ring $Q=\Znil[S]$ is
isomorphic to the category of pairs $(G,\alpha)$, where $G$ is a
group of nilpotence class two and $\alpha:S\to \Hom(G,G)$ is an
action of $S$ on $G$ via group homomorphisms.
\subsection{Internal groupoids and crossed objects}
We describe now internal groupoids in the category of square groups.
Actually results obtained in this section are very particular case
of much more general results of Gogi Janelidze \cite{gogia}.

Let $A$ and $G$ be square groups. An \emph{action of $G$ on $A$} is
a homomorphism of abelian groups $\xi:A\ad\t G\ad\to A_{ee}$.

In particular we have the action of $A$ on $A$ given by $(-,-)_H$,
which is called the \emph{adjoint action of $A$ on itself}.

Let $\xi$ be an action of $G$ on $A$. The \emph{semi-direct} product
of $G$ and $A$ denoted $G\rtimes A$ is a square group defined as
follows. As a set $(G\rtimes A)_e$ is the cartesian product $G_e\x
A_e$ while the group structure is given by
$$(g,x)+(h,y)=(g+h,x+y+P\xi(g,h)).$$
Moreover one puts
$$(G\rtimes A)_{ee}=G_{ee}\oplus A_{e},$$
$$P(u,a)=(Pu,Pa),$$
$$H(g,x)=(H(g),H(x)-\xi(x,g)).$$
One easily sees that
$$[(g,x),(h,y)]=([g,h],[x,y]+P\xi(x,h)-P\xi(y,g))$$
and
$$((g,x)\mid (h,y))_H=((g\mid)_H,(x\mid y)_H+HP\xi(x,h)-P\xi(x,h)-P\xi(y,g)).$$
Based on these identities one easily checks that $(G\rtimes A)_e$ is
really a square group and one has the following split short exact
sequence of square groups
$$0\to A\to A\rtimes G\to G\to 0$$
with obvious maps. Conversely, let
$$0\to A\xto{i} B\xto{p} G\to 0$$
be a short exact sequence. Assume $j:G\to B$ is  a morphism of
square groups with $pj=\Id_G$. Then
$$\xi(x,g):=(i_e(x)\mid j_e(g))_H$$
defines an action of $G$ on $A$ and the maps
$f_e(g,x)=j_e(g)+i_e(x)$ and $f_{ee}(g,x)=j_{ee}(g)+i_{ee}(x)$
define an isomorphism $f=(f_e,f_{ee}):G\rtimes A\to B$ of square
groups.

A \emph{crossed square group} is a morphism of square groups
$\d:A\to G$ together with an action of $G$ on $A$ such that $\d$ is
compatible with the action of $G$, where $G$ acts on itself via the
adjoint action and the action of $A$ on $A$ given via $\d$ coincides
with the adjoint action of $A$ on itself. In other words a
homomorphism $\xi:A\ad\t G\ad\to A_{ee}$ of abelian groups is given
and the following identity holds
$$\d_{ee}\xi(x,g)=(\d_e(x),g)_H,$$
$$\xi(x,\d_e(y))=(x,y)_H.$$
We denote by $\CSG$ the category of crossed square groups. The
following is a specialization of the main result of \cite{gogia}.

\begin{Le}
Any internal category in the category of square groups is an
internal groupoid. Thus $\Cat(\SG)=\Gpd(\SG)$ and there is an
 equivalence of categories
$$\Gpd(\SG)\cong \CSG.$$
\end{Le}

\begin{proof} The first fact is a general property of so
called Maltsev categories \cite{gogia}. The second part can be
proved by modifying the argument of Loday in \cite{jllcat} based on
our description of split short exact sequences. Alternatively one
can check directly that the above definition is a specialization of
the general notion of Janelidze and use the main result of
\cite{gogia} on relationship between internal groupoids and crossed
objects in so called semi-abelian categories. The checking is an
easy exercise because of the explicit description of the coproduct
in the category $\SG$ of square groups given in \cite{square}.
\end{proof}

Since the functors $(-)_e:\SG\to \Gr$ and $(-)_{ee}:\SG\to \Ab$
preserve limits any internal groupoid $X$ in $\SG$ gives rise to two
internal groupoids $X_e$ and $X_{ee}$ in the category of groups and
abelian groups respectively.

An internal groupoid $X\in\Gpd(\SG)$ is called ee-antidiscrete
provided $X_{ee}$ is antidiscrete. We denote by $\Gpd_{adee}(\SG)$
the category of ee-antidiscrete internal groupoids in the category
of square groups.

\begin{Le}\label{lodjan}
The equivalence $\Gpd(\SG)\cong {\sf cross}(\SG)$ restricts to an
equivalence of categories $$\Gpd_{adee}(\SG)\cong \QPM.$$
\end{Le}

\begin{proof}
We have to show that qpms are exactly  crossed square groups
$\d:A\to G$ for which $\d_{ee}:A_{ee}\to G_{ee}$ is the identity
map. But this is clear, because after identification of $A_{ee}$ and
$G_{ee}$ via $\d_e$, the action $\xi$ of $G$ on $A$ becomes
redundant, $\xi(x,g)=(\d_e(x)\mid x)_H$.
\end{proof}

Let $\Gpd_{adee}(\SR)$ denote the category of ee-antidiscrete
groupoid objects in the category of square rings.  Lemma
\ref{lodjan} implies the following result.

\begin{Le}\label{lodgog} There is an equivalence of categories
$$\Gpd_{adee}(\SR)\cong \CS.$$
\end{Le}
\section{Proof of the main result}
\subsection{Relative Mac~Lane cohomology}
Let $R$ be a ring and let $M$ be a bimodule over $R$. Assume also that a
  surjective morphism $p:Q\to R$ is given  from a square
ring $Q$ to $R$. We denote by $\xext(R,Q;M)^{\SG,\square}$ the
subcategory of the category $\xext(R,M)^{\SG,\square}$ whose objects
are crossed square ring extensions  of the form
$$\xymatrix{
0\ar[r]& M\ar[r]& C_{(1)}\ar[r] ^{(\d,\Id)} &
C_{(0)}\ar[r]^p&R\ar[r]&0}$$ with
$C_{(0)}=Q$. Morphism are such morphisms of crossed square ring
extensions which are identity on $Q$
$$\xymatrix{
0\ar[r]& M\ar[r]\ar[d]_\Id& C_{(1)}\ar[d]^f\ar[r] ^{(\d,\Id)} &Q
\ar[d]^\Id\ar[r]^p&R\ar[d]^\id
\ar[r]&0\\
0\ar[r]& M\ar[r]& C_{(1)}'\ar[r] ^{(\d',\Id)} &Q
\ar[r]^p&R\ar[r]&0.}$$ Then the category
$\xext(R,Q;M)^{\SG,\square}$ is a groupoid.

Quite similarly, for a given surjective morphism $p:Q\to R$  from a
quadratic ring $Q$ to $R$, we denote by $\xext(R,Q;M)^{\SG,\tl}$ the
subcategory of the category $\xext(R,M)^{\SG,\tl}$ whose objects are
crossed square ring extensions of the form
$$
\xymatrix{
0\ar[r]& M\ar[r]& C_{(1)}\ar[r] ^{(\d,\Id)} &C_{(0)} \ar[r]^p&R\ar[r]&0
}
$$
with $C_{(0)}=Q$. Then the category  $\xext(R,Q;M)^{\SG,\tl}$ is a
groupoid.

\begin{Le}\label{monreliso} Let $R$ be a ring and let $L$ be a monoid
and let $p:\Znil[L]\to R$ be a surjective morphism of quadratic
rings (and hence also a surjective morphism of square rings). Then
for any $R$-bimodule $M$ the functor $\xext(R,M)^{\SG,\tl}\to
\xext(R,M)^{\SG,\square}$ yields an equivalence of categories
$$\xext(R,\Znil[L];M)^{\SG,\tl}\xto\simeq\xext(R,\Znil[L]
;M)^{\SG,\square}.
$$
\end{Le}

\begin{proof} It is straightforward to check that
the conditions posed on $C_{1}$ and $\d$ in the definition of
quadratic pair algebra and in the definition of crossed square ring
are the same provided $C_{(0)}=\Znil[L]$.
\end{proof}

Let us turn back to an epimorphism $Q\to R$ for a square ring $R$.
The set of connected components of $\xext(R,Q;M)^{\SG,\square}$ has
the following cohomological description. In order to give the
precise statement we first extend the definition of the Mac~Lane
cohomology to square rings and then we introduce the relative
cohomology groups.

Let $Q$ be a square ring, then $Q\ad$ is a ring, which we denote by
$R$. There is an obvious functor $$q:\modr Q\to \modr R$$ which is
identity on objects and on morphisms it is given by
$$q((f_i^k,f_{ij}^k)):=(\bar{f}_i^k)$$ For any bimodule $M$ over the
ring $R$ we let $D_M$ be the bifunctor on $\modr R$ given by
$$(X,Y)\mapsto \Hom_R(X,Y\t _R M)$$
By abuse of notation we will denote by $D_M$ also the induced bifunctor on
$\modr Q$. Then we put
$$
\HML^*(Q,M):=H^*(\modr Q,D_M)
$$
Thanks to Section \ref{relcoh} we recover for usual rings the classical Mac
Lane cohomology. Using the relative cohomology of small categories defined in
Section \ref{relcoh} one can also define the relative Mac Lane cohomology groups
$\HML^*(R,Q;M)$ to be $H^*(\modr R,\modr Q;D_M).$ Thus one has the following
long exact sequence
\begin{multline*}
0\to \HML^0(R;M)\to \HML^0(Q;M)\to \HML^1(R,Q;M)\to\cdots\\
\to \HML^n(R;M)\to \HML^n(Q;M)\to \HML^{n+1}(R,Q;M)\to\cdots.
\end{multline*}

The proof of the isomorphisms
in Theorem \ref{xext} is based on a computation given in
Appendix and on the following result.
\begin{The}\label{fardobitisq}
Let $p:Q\to R$ be a surjective morphism from a square ring $Q$ to a ring
$R$.
Then
$$\pi_0(\xext(R,Q;M)^{\SG,\square})\approx \HML^3(R,Q;M).
$$
\end{The}
\begin{proof}
Let ${\Tracks}(\modr R, \modr Q;D_M)$ denote the category of such
abelian track categories $\ta$ that the corresponding homotopy
category $\ta_\ho$ is $\modr R$, underlying category $\ta_0$ is
$\modr Q$ and the corresponding natural system is given by the
bifunctor $D_M$. We now construct the functor
$$
\chi:\xext(R,Q;M)^{\SG,\square}\to {\Tracks}(\modr R, \modr Q;D_M)
$$
as follows. Let
$$
\xymatrix{
0\ar[r]& M\ar[r]& \tilde{Q}\ar[r] ^{(w,\Id)} &Q \ar[r]^p&R\ar[r]&0,
}
$$
be a crossed square ring extension. The underlying category of the
track category $\chi(\omega)$ is $\modr Q$. If ${\bf
f}=(f^k_i,f^k_{ij})$ and ${\bf g}=(g^k_i,g^k_{ij})$ are morphisms
$y\to x$, $x,y\in \N$ in $\modr Q$, then a track $\bf f\then g$ is a
collection $(h^k_i)$ of elements in $\tilde{Q}_e$ such that
$\d(h^k_i)=f^k_i-g^k_i$ for all $1\le i\le x$ and $1\le k\le y$. Now
the result follows from the fact that $\chi$ is an isomorphism of
categories. The inverse of $\chi$ is given  as follows. Let $\ta$ be
a track category such that $\ta_\ho=\modr R$ and $\ta_0=\modr Q$. By
Lemma \ref{lec} $\ta$ is an internal groupoid in the category of
quadratic theories. By applying the functor $\sf cro$ one obtains an
internal groupoid in the category of square rings. Moreover, the
proof of Lemma \ref{lec} shows that this groupoid is
ee-antidiscrete, therefore by Lemma \ref{lodgog} it defines an
object in $\xext(R,Q;M)^{\SG,\square}$, which is the value of the
inverse of $\chi$.
\end{proof}

\subsection{A pullback construction}\label{pbc}
We now give a construction in the category of crossed square ring
extensions which is needed in the proof of Theorem \ref{xext}.
Let
$$\xymatrix{
0\ar[r]& M\ar[r]& C_{(1)}\ar[r] ^{(\d,\Id)} &C_{(0)}\ar[r]^p&R\ar[r]&0,}
$$
be a crossed square ring extension and let $f:Q_{(0)}\to C_{(0)}$ be a
morphism of square rings, such that $p\circ f_e:Q_{0}\to R$ is
surjective. Based on this data we construct  the following crossed
square ring
$$
\xymatrix{& Q_{ee}\ar[dl]_{P^Q}\\
Q_{1}\ar[r]_{\d^Q} &Q_{0}\ar[u]_{H^{Q_0}}}
$$
where the group $Q_{1}$ is defined by the pullback diagram
$$\xymatrix{Q_{1}\ar[r]^{\d^Q}\ar[d]_{g_e}&Q_{0}\ar[d]^{f_e}\\
C_{1}\ar[r]_\d& C_{0}}$$ and $P^Q=(P^{C}\circ f_{ee},
P^{Q_0}):Q_{ee}\to Q_{1}$. Then one has the following crossed square
ring extension
$$\xymatrix{
0\ar[r]& M\ar[r]& Q_{(1)}\ar[r] ^{\d^Q} &Q_{(0)} \ar[r]&R\ar[r]&0.}$$
One easily sees that
$$\xymatrix{
0\ar[r]& M\ar[r]\ar[d]_\Id& Q_{(1)}\ar[d]^g\ar[r] ^{\d^Q} &Q_{(0)}
\ar[d]^f\ar[r]^{pf}&R\ar[d]^\id
\ar[r]&0\\
0\ar[r]& M\ar[r]& C_{(1)}\ar[r] ^{\d} &C_{(0)} \ar[r]^p&R\ar[r]&0}$$ is a morphism
of crossed square ring extensions.

We call this construction the pullback construction and write
$f^*\d$ instead of $(\d^Q)$. Assume now that  $(\d)$ is linearly
generated and the composite ${\sf L}(Q_{(0)})\to {\sf L}(C_{(0)})\to R$ is
surjective, then one easily sees that $(f^*\d)$ is also linearly
generated.

Of course a similar constructions works for  quadratic pair algebras.

\subsection{Proof of Theorem \ref{xext}}\label{proof261} Let
$$
\xymatrix{
0\ar[r]& M\ar[r]& \tilde{Q}\ar[r] ^{(w,\Id)} &Q \ar[r]^q&R\ar[r]&0,}
$$
be an object of $\xext_L(R,M)^{\SG,\square}$. For simplicity we
denote this object by $(w)$. Then it can be also considered as an
object of $\xext(R,Q;M)^{\SG,\square}$ and therefore $(w)$ defines
an element in $\HML^3(R,Q;M)$ thanks to Theorem \ref{fardobitisq}.
Then the boundary homomorphism gives an element in $\HML^3(R,M)$. In
this way we get a map
$$
\zeta: \pi_0(\xext_L(R,M)^{\SG,\square})\to \HML^3(R;M).
$$
Composing it with $\pi_0(\xext_L(R,M)^{\SG,\tl}) \to
\pi_0(\xext_L(R,M)^{\SG,\square})$ we obtain the map
$$
\zeta': \pi_0(\xext_L(R,M)^{\SG,\tl}) \to \HML^3(R;M).
$$
We have to show that these maps are bijections. Take an $a\in
\HML^3(R;M)$. Take any  surjective homomorphism $L\to R$ from a free
monoid $L$ to the multiplicative monoid of the ring $R$. It yields a
surjective morphism $r:\Znil[L]\to R$. Here $\Znil[L]$ can be
considered as a square ring as well as a quadratic ring. Since
$\HML^i(\Znil[L];D_M)=0$ for $i=2,3$ (see Theorem \ref{nilnul} in
Appendix), we have an isomorphism
$$
\partial:\HML^3(R,\Znil[L];M)\cong\HML^3(R;M).
$$
Let $b=\partial\1(a)\in\HML^3(R,\Znil[L];M)$ be the element
corresponding to $a$. Thanks to Theorem \ref{fardobitisq} the
element $b$ defines a crossed square ring extension of $R$ by $M$
$$
\xymatrix{
0\ar[r]& M\ar[r]& \tilde{Q}\ar[r]^-{(v,\Id)}&\Znil[L]
\ar[r]&R\ar[r]&0}
$$
which is also linearly generated by construction and therefore is an
object of $\xext_L(R,M)^{\SG,\square}$. By Lemma \ref{monreliso} it
can be considered also as a quadratic pair algebra extension. Hence
$\zeta$ and $\zeta'$ are surjections. It remains to show that
$\zeta$ and $\zeta'$ are injections as well. Suppose
$\zeta(w)=\zeta(w')$ (resp. $\zeta'(w)=\zeta'(w')$). We have to show
that $(w)$ and $(w')$ are in the same connected component. Let ${\sf
L}(Q)$ be the monoid of linear elements in $Q$. Via $q$ it maps to
the multiplicative submonoid $q({\sf L}(Q))$ of $R$. Take any
surjective homomorphism of monoids $F\to q({\sf L}(Q))$ with $F$ a
free monoid. It has a lifting to a monoid homomorphism $F\to {\sf
L}(Q)$, which yields a square (resp. quadratic) ring homomorphism
$t:\Znil[F]\to Q$. The homomorphism $t$ satisfies all conditions on
$f$ in Section \ref{pbc} and hence yields a morphism of crossed
square ring extensions (resp. quadratic pair algebra extensions)
$t^*(w)\to w$. Thus without loss of generality we can assume that
$(w)$ and $(w')$ are chosen in such a way that $Q=\Znil[F]$ and
$Q'=\Znil[F']$.  Let $L$ and $r$ be the same as above (see the proof
of surjectivity of $\zeta$). Since $L\to R$ is surjective,
$q(F)\subset R$ and $F$ is free, there exists a morphism of monoids
$F\to L$ such that for the induced morphism $k:Q=\Znil[F]\to
\Znil[L]$ one has $q=r\circ k$. Thus one has the following
commutative diagram
$$\xymatrix{\HML^3(R,\Znil[L];M)\ar[r]\ar[d]^{k^*}& \HML^3(R,M)\\
 \HML^3(R,Q;M)\ar[ur]}$$
Since both morphisms in the diagram  with target $HML^3(R,M)$ are
isomorphisms, it follows that $k^*:\HML^3(R,\Znil[L];M)\to
\HML^3(R,Q;M)$ is also an isomorphism. Considering an extension
corresponding to ${k^*}^{-1}(w)$ one sees that there exists a
morphism of square ring extensions (resp. quadratic pair algebra
extensions)
$$\xymatrix{0\ar[r]& M\ar[r]\ar[d]^{\Id}& \tilde{Q}\ar[d]\ar[r]^-{(w,\Id)}
&Q\ar[d]^k
\ar[r]^q&R\ar[r]\ar[d]^\Id&0\\
0\ar[r]& M\ar[r]& \bar{Q}\ar[r]^-{(\bar{w},\Id)} &\Znil[L]
\ar[r]&R\ar[r]&0}$$ In a similar manner we find a morphism of square
ring extensions (resp. quadratic pair algebra extensions)
$$\xymatrix{0\ar[r]& M\ar[r]\ar[d]^{\Id}& \tilde{Q'}\ar[d]\ar[r]^-{(w',\Id)}
&Q\ar[d]^k
\ar[r]^{q'}&R\ar[r]\ar[d]^\Id&0\\
0\ar[r]& M\ar[r]& \bar{Q'}\ar[r]^-{(\bar{w'},\Id)} &\Znil[L]
\ar[r]&R\ar[r]&0}$$ Since the square ring extensions (resp.
quadratic pair algebra extensions) $(\bar{w})$ and $(\bar{w'})$ lie
in the same groupoid $\xext(R,\Znil[L] ;M)^{\SG,\square}$ and their
classes in $ \HML^3(R,\Znil[L];M)$ are the same, it follows that
they are isomorphic in the groupoid $\xext(R,\Znil[L]
;M)^{\SG,\square}$. Therefore we have the following diagram in $
\xext(R,M)^{\SG,\square}$ (resp. $\xext(R,M)^{\SG,\tl} $):
$$
(w')\ot  (\bar{w}')\ \cong (\bar{w})\to (w),
$$
hence the result.

\section{Application to theory of
2-categories}\label{obs2cat}
In this section we introduce the notion of 2-additive track
category, which is the 2-categorical analogue of
additive category and we prove a strengthening theorem for such
2-additive track categories.

\subsection{Abelian track categories}
A track category is \emph{abelian} if for any map $f:X\to Y$, the
group $\Aut(f)$ of tracks from $f$ to itself is abelian.  Any track
category which fits in a linear track extension is abelian. Converse
is also true: any abelian  track category defines a natural system
$D=D_{\ta}$ on $\ta_\ho$ and a linear track extension
$$
0\to D_\ta\to \ta_1 \rightrightarrows \ta_0\to \ta_\ho \to 0.
$$
The natural system $D_\ta$ and the linear track extension are unique
up to isomorphism (see \cite{BJ}).

\subsection{Track theories}
A \emph{coproduct} $A\vee B$ in a track category $\ta$ is an object
$A\vee B$ equipped with 1-morphisms $i_1:A\to A\vee B$, $i_2:B\to
A\vee B$ such that the induced functor
$$
(i_1^*,i_2^*):\hog{A\vee B,X}\to \hog{A,X}\x \hog{B,X}
$$
is an equivalence of groupoids for all objects $X\in\ta$. The
coproduct is \emph{strong} if the functor $(i_1^*,i_2^*)$ is an
isomorphism of groupoids.  By duality we have also notion of
\emph{product} and  \emph{strong product}.  A \emph{zero object} in
a track category $\ta$ is an object $0$ such that the categories
$\hog{0,X}$ and $\hog{X,0}$ are equivalent to the trivial groupoid
for all $X\in\ta$.
A \emph{strong zero object} in a track category $\ta$ is an object $0$
such that all categories $\hog{0,X}$ and $\hog{X,0}$ are trivial
groupoids.

A \emph{track theory} (resp. \emph{strong track theory}) is a small track
category $\ta$ possessing finite coproducts (resp. strong coproducts).
Morphisms of track theories are enriched functors which are compatible with
coproducts. An equivalence of track theories is a track theory morphism which is
a weak equivalence and two track theories are called \emph{equivalent} if they
are made so by the smallest equivalence relation generated by these.
The following is a particular case of a general result of Power
\cite{Pow}. For a cohomological proof we refer to  \cite{streng}.

\begin{The}\label{power}
Any abelian track theory is equivalent to a strong one.
\end{The}

If $\ta$ is an abelian track theory, then the corresponding category
$\ta_\ho$ is an algebraic theory and the natural system $D_\ta$ is
cartesian. Conversely, if $\ta$ is an abelian track category such
that $\ta_\ho$ is an algebraic theory and $D_\ta$ is cartesian, then
$\ta$ is an abelian track theory.
Moreover an abelian track theory
is strong if and only if $\ta_0$ is an algebraic theory and the
canonical functor $\ta_0\to\ta_\ho$ is a morphism of algebraic
theories.

\subsection{2-Additive track categories}
Now we introduce an analogue of additive categories in the 2-world.
Let $\ta$ be a track theory with zero object. Then for any objects
$A$ and $B$ of $\ta$, there is a map $p_1: A\vee B\to A$ and tracks
$p_1i_1\then \id_A$, $p_1i_2\then 0$. Similarly for $p_2:A\vee B\to B$.
A \emph{2-additive track category} is an abelian track
theory with strong zero object, such that the following conditions hold

 i) for any two objects $A$
and $B$ the coproduct $A\vee B$ is also a product via $p_1:A\vee B\to A$
and $p_2:A\vee B\to  B$

ii)  for any morphism $f:A\to B$ there exists a morphism $g:A\to B$ and
a track $hd\then 0$, where $d:A\to A\vee A$ and $h:A\vee A\to B$ are
morphisms with tracks $hi_1\then f$, $hi_2\then g$, $p_1d\then \id_A$,
$p_2d\then \id_A$.

It is clear that the homotopy category
$\ta_\ho$ of a 2-additive track theory is an additive
category. The following is a direct consequence
of \cite{BJ} and \cite{beitrage}.

\begin{Le}
Let $\ta$ be an abelian track category. Then $\ta$ is a 2-additive
track category iff $\ta_\ho$ is an additive category and the
corresponding natural system $D_\ta$ is a biadditive bifunctor.
\end{Le}

It follows that a 2-additive track category determines a triple
$(\ta_\ho,D_\ta,\Ch(\ta)\in H^3(\ta_\ho;D_\ta))$. Conversely for an
additive category $\C$, a biadditive bifunctor  $D:\C\op\x\C\to\Ab$
and an element $a\in H^3(\C;D)$ there exists a 2-additive track
category unique up to equivalence such that $\ta_\ho=\C$, $D_\ta=D$
and $\Ch(\ta)=a$.

\subsection{Strongly and very strongly 2-additive track theories}
As we said Theorem \ref{power} asserts that any track theory is
equivalent to one with strong coproducts. In particular, any
2-additive track category is equivalent to one which possesses
strong products. Since the dual of an additive track category is
still a track theory, we see that it is also equivalent to one which
possesses strong coproducts. Can we always get strong products and
coproducts simultaneously? In other words, is every 2-additive track
category $\ta$ equivalent to a very strongly 2-additive track
theory? Here a 2-additive track category is called \emph{very strongly
2-additive} if it admits a strong zero object $0$, strong finite
coproducts and for any two object $A$ and $B$ the strong coproduct
$A\vee B$ is also a strong coproduct via $p_1:A\vee B\to A$ and
$p_2:A\vee B\to B$. The answer is given by the following result, which also
shows that the number 2 plays an important r\^ole
in the theory of 2-categories.

\begin{The} \label{verystongstrengthening}
Let $\ta$ be a small 2-additive track category with homotopy
category $\C=\ta_\ho$ and canonical bifunctor $D=D_{\ta}$. Let $_2D$
be the two-torsion part of $D$. Then there is a well-defined element
$\nu(\ta)\in H^0(\C;{}_2D)$, which is nontrivial in general and such
that $\nu(\ta)=0$ iff $\ta$ is equivalent to a very strongly
2-additive track theory. The class $\nu(\ta)$ for example is zero
provided homs of the additive category $\C$ are modules either over
$\Z[\frac12]$ or over $\F_2$ (the field with two elements).
\end{The}

\begin{proof} First one observes that a 2-additive track category
$\ta$ is very strongly 2-additive iff the category $\ta_0$ is
additive. For simplicity we restrict ourself to the case of single
sorted theries. Then $\ta_\ho=\modr R$ for a ring $R$. In this case
one has an isomorphism $D(X,Y)\cong\Hom(X,M\t_R Y)$ natural in $X$
and $Y$, where $M=D(R,R)$ is a bimodule over $R$. We claim that up
to equivalence single sorted very strongly 2-additive track
categories $\ta$ with fixed $\ta_\ho$ and $D_\ta$ are in bijection
with $\Sh^3(R;M)$. Indeed, if
$$
0\to M\to C_1\xto\d C_0\to R\to 0
$$
is an object of ${\sf Cross}(R,M)$, then we have the following very
strongly 2-additive track category $\ta$. Objects of $\ta$ are the
same as the objects of $\modr R$, i.~e. natural numbers. For any
natural numbers $n$ and $m$ the maps from $n$ to $m$ (which is the
same as objects of the groupoid $\ta(n,m)$)  are $m\times
n$-matrices with coefficients in $C_0$. For $f,g\in Mat_{m\times
n}(S)$ the set of tracks $f\to g$ is given by
$$
\Hom_{\ta(n,m)}(f,g)=\left\{h \in \mathrm{Mat}_{m\times n}(C_1)\
\mid\ \d (h)=f-g\right\}.
$$
Composition of 1-arrows is given by the usual multiplication of
matrices, while composition of tracks is given by the addition of
matrices. One easily checks that in this way one really  obtains a
very strongly 2-additive track theory $\ta(\d)$. It is clear that
$\ta_\ho=\modr R$, where $R={\sf Coker}(\d)$ and the bifunctor
corresponding to $\ta$ is $D=\Hom(-,M\t_R-)$. Conversely, assume
$\ta$ is a single sorted very strongly 2-additive track category
with $\ta_\ho=\modr R$ and $D_\ta=\Hom(-,M\t_R-)$. Then $\ta_0$ is a
single sorted additive category and therefore it is equivalent to
$\modr S$, where $S={\sf End}_{\ta_0}(1)$. Restriction of the
quotient functor $\ta\to \ta_0$ yields a homomorphism of rings $S\to
R$. One defines $X$ to be the set of pairs $(h,x)$, where $x\in
\Hom_{\ta_0}(1.1)$ and $h:x\then 0$ is a track in the groupoid
$\ta(1,1)$. Moreover we put $\partial=\partial_{\ta}(h,x)=x$. Then
$X$ carries a structure of a bimodule over $S$, and
$$
0\to M\to X\buildrel\partial\over\to S\to R\to 0
$$
is a crossed extension and the claim follows from isomorphism
(\ref{shuxt}). As we said, up to equivalence single sorted 2-additive track
categories $\ta$ with fixed $\ta_\ho$ and $D_\ta$ are in bijection
with $\HML^3(R,M)$. Therefore the exact sequence $0\to \Sh^3(R,M)\to
\HML^3(R,M)\xto{\nu} H^0(R, \,_2M)$ together with Proposition 9.1.1 of
\cite{shukla} implies the result.
\end{proof}

{\bf Remark.} One can describe the function $\nu$ in Theorem
\ref{verystongstrengthening} as follows. Let $\ta$ be a 2-additive
track theory. Let $\vee$ denote the weak coproduct in $\ta$ and let
$0$ be the weak zero object. For objects $X,Y$ one has therefore
``inclusions'' $i_1:X \to X\vee Y$ and $i_2:Y \to X\vee Y$. Since
$X\vee Y$ is also a weak product of $X$ and $Y$ in $\ta$ it follows
that one has also projection maps $p_1:X\vee Y\to X$ and $p_2:X\vee
Y\to Y$.
For each $X$ we choose maps $i_X:X\to X\vee X$ and $t_X:X\vee Y\to
Y\vee X$ in such a way that classes of $i_X$ and $t_X$ in $\ta_\ho$
are the codiagonal and twisting maps in the additive category
$\ta_\ho$. It follows that there is a unique track
$$
\aa_X:i_X\then t\circ i_X
$$
such that $p_{i*}(\aa_X)=0$ for $i=1,2$. Now, let $(1,1):X\vee X\to
X$ be a map which lifts the codiagonal map in $\ta_\ho$. Then
$(1,1)_*\aa_X$ is a track $\id_X\to\id_X$ and therefore it differs
from the trivial track by an element $\nu(X)\in D(X,X)$. One can
prove that the assignment $X\mapsto\nu(X)$ is the expected one.

\subsection{Strongly additive track categories}

A 2-additive track category $\ta$ is called \emph{strongly
2-additive} if $\ta_0$ is quadratic.

\begin{The}\label{stronglyadditive}
Any 2-additive track theory is equivalent to a strongly 2-additive one.
\end{The}

\begin{proof}
We continue to restrict ourselves to the single sorted case. In this
case $\ta_\ho$ is the category $\modr R$ for a ring $R$ and
$D=\Hom_R(-,(-)\t_R M)$ for an $R$-bimodule $M$. Thus ${\sf
Ch}((\ta))\in \HML^3(R;M)$ and therefore it belongs to
$\HML^3(R,Q;M)$ for a square ring $Q$ thanks to the proof of Theorem
\ref{xext} given in Section \ref{proof261}. Thus the element ${\sf
Ch}((\ta))$ has a realization via track category $\ta'$ such that
$(\ta')_0=\modr Q$ and we are done.
\end{proof}


\appendix

\section{Cohomology of free monoid square rings}

\centerline{\sc T. Pirashvili}

\

\
Here we prove the following result.

\begin{The} \label{nilnul}
Let $L$ be a free monoid and let $Q=\Znil[L]$ be the corresponding
monoid square ring. Then for any $R$-$R$-bimodule $B$ one has
$$
\HML^2( Q,B)=0= \HML^3( Q,B).
$$
\end{The}

Proof of Theorem \ref{nilnul} is given in Section \ref{proofnil}.
The argument is a modification of the one given in \cite{Nil}.

\subsection{Auxiliary results}
For a ring $R$ we denote by ${\bf F}(R)$ or simply by $\bf F$ the
category of all covariant functors from the category $\modf R$ of
finitely generated free right $R$-modules to the category $\Mod R$
of all right $R$-modules. It is well known \cite{JP} that
$$
\HML^*(R,B)\cong \Ext^*_{\bf F}(\Id, (-)\t_RB)
$$

We need the following result, which is an easy consequence of
Theorem 9.2.1 \cite{shukla} and the fact that ${\sf SH}^i(R,-)=0$
for all $i\ge 2$, provided $R$ is a free ring.

\begin{Le}\label{freering} Let $R$ be a free ring and let $B$ be an
$R$-$R$-bimodule. Then one has $\HML^2(R, B)=0$
 and $\HML^3(R, B)\cong {\sf H}^0 (R,\, _2
B)$,  where ${\sf H}^*(R,-)$ denotes the Hochschild cohomology of $R$.
\end{Le}

 We also need the following vanishing result.

\begin{Le}\cite{add}\label{vanishing}
Let $R$ be a ring and let $$T:(\modf R)\x (\modf R)\to \Mod R$$ be a
bifunctor, which is covariant in both variables and
$T(0,X)=0=T(X,0)$ for all $X\in \modf R$. Then for any additive
functor $F:\modf R\to \Mod R$ one has
$$
\Ext^*_{\bf F}(F, T^d)=0=\Ext^*_{\bf F}(T^d,F),
$$
where $T^d(X)=T(X,X)$.
\end{Le}

In the following we need the simplicial derived functors of the
functor $(-)\ad:Q$-$Mod\to R$-$Mod$, which are denoted by
$$
\Tor^Q_*(-,R):\Mod Q\to \Mod R.
$$
We recall the definition of these functors. According to
\cite{quillen} the category of simplicial objects in the category of
right $Q$-modules has a closed model category structure where a
morphism $f:X_*\to Y_*$ of simplicial objects is a weak equivalence
(resp. fibration) when it is so in the category of simplicial sets.
Let $M$ be a right $Q$-module and let $X_*$ be a cofibrant
replacement of $M$.  By \cite{quillen} one can assume that each
$X_n$, $n\ge 0$ is a free right $Q$-module. We also have
$\pi_iX_*=0$ for $i>0$ and $\pi_0X_*=M$. Now one puts
$$\Tor^Q_*(M,R):=\pi_*(X_*\ad).$$
It is well known that these are well-defined functors. Since $\Mod
R\subset \Mod Q$, one can consider also the restriction of
$\Tor^Q_*(-,R)$ to $\Mod(R)$ (see Proposition \ref{ssnil} below).

\begin{Pro}\label{ssnil}
For any square ring $Q$ and for any $R$-$R$-bimodule $B$, one has the following
spectral sequence
$$E^2_{pq}=\Ext^p_{\bf F}(\Tor_q^Q(-,R),F)\Longrightarrow \HML^{p+q}( Q,B), $$
where $R=Q\ad$, $F(-)=(-)\tp _R B$.
\end{Pro}

{\it Proof}. Proposition follows immediately  from the spectral sequence
(8.2.2) and Lemma 8.3.1 of \cite{Nil}.

\subsection{Computation of $\Tor^Q$} In this section we give a computation of
$\Tor$-groups involved in Proposition \ref{ssnil}. It is based on
Lemma \ref{qucfo} below, which is the specialization of the exact
sequence (4.1) of \cite{qucf}. Let us recall that Eilenberg and
Mac~Lane \cite{EM} defined the quadratic functor
$$\Omega:\Ab\to \Ab$$ such that it commutes with filtered colimits,
$$\Omega(A\oplus B)=\Omega(A)\oplus  \Omega(B)\oplus \Tor(A,B)$$
and moreover $$\Omega(\Z)=0, \ \ \Omega(\Z/n\Z)=\Z/n\Z.$$

\begin{Le}\cite{qucf} \label{qucfo} Let $X_*$ be a simplicial abelian group, which is degreewise free and
has homotopy groups $\pi_i=\pi_i(X_*)$. Then one has
$$\pi_0(\La ^2X_*)=\La ^2 (\pi_0)$$
$$0\to \pi_1\t \pi_0\to \pi_1(\La^2(X_*))\to \Omega(\La ^2\pi_0)\to 0$$
$$0\to \pi_2\t \pi_0 \oplus \Gamma \pi_1
\to \pi_2(\La ^2 (X_*))\to \Tor( \pi_1,\pi_0)\to 0$$
\end{Le}

Let us recall that if $G$ is a free class two nilpotent group, then
one has the following short exact sequence
$$0\to \La^2(G\ab)\to G\to G\ab\to 0,$$
where the first nontrivial map is induced by $(x,y)\mapsto
-x-y+x+y$. Assume now that $L$ is a monoid and  $Q=\Znil[L]$ is the
corresponding monoid square ring. As we already mentioned a right
$Q$-module is the same as a nilpotent group of class two together
with an action of $L$ via group homomorphisms. It follows that if
$X$ is a free right $Q$-module, then $X$ is  also  free as a
nilpotent group of class two. Furthermore, $X\ad$ in this case is
simply $X\ab$, thus we have the following Lemma.

\begin{Le}\label{sesla} Let $L$ be a monoid and let $Q=\Znil[L]$ be
the monoid square ring.  Then, for any free right $Q$-module $X$, one has the
following short exact sequence
$$0\to \La^2(X\ad)\to X\to X\ad\to 0,$$
in the category of modules over the ring $R=Q\ad=\Z[L]$, where the
first nontrivial map is induced by $(x,y)\mapsto -x-y+x+y$, and
$\La^2(X\ad)$ is an $R$-module via the diagonal action of $L$.
\end{Le}

We would like to use these results in the following situation.
\begin{Pro}\label{qtor} Let $L$ be a monoid and let $Q=\Znil[L]$ be
the monoid square ring. Then, for any free right module $M$ over the
ring $R=Q\ad$, one has the following natural isomorphisms
$$\Tor_0^Q(M,R)\cong M$$
$$\Tor_1^Q(M,R)\cong \La^2(M)$$
$$\Tor_2^Q(M,R)\cong M\tp \La^2(M)$$
$$\Tor_3^Q(M,R)\cong (M\tp M\tp \La^2(M))\oplus (\Gamma (\La ^2(M)))$$

\end{Pro}

{\it Proof}. Let $M$ be a free $R$-module. Let us take a free simplicial
resolution $Y_*$  of $M$ in the category of $Q$-modules. Thanks to Lemma
\ref{sesla} one has an exact sequence
$$0\to \La ^2X_*\to Y_*\to X_*\to 0,$$
where $X_*=Y_*\ad$. Since $\pi_iY_*=0$ for $i>0$ and
 $\pi_0Y_*=M$ we have $\pi_0X_*=M$ and
$\pi_{i+1}X_*= \pi_i\La ^2(X_*)$. Since  $M$ is a free abelian
group, one can use Lemma \ref{qucfo} to get
$$\pi_1(X_*)\cong \La ^2(M), \ \ \  \pi_2 (X_*)= M\t \La ^2 (M)$$
$$\pi_3(X_*)\cong  \La ^2(M)\t M^{\t 2}\oplus \Gamma (\La ^2M)$$
Comparing with definition of simplicial derived functors we obtain  the
expected result.

 \subsection{Universal quadratic functors}
Let $A$ be an abelian group. We set
$$P(A)=I(A)/I^3(A),$$
where $I(A)$ is the augmentation ideal of the group algebra of $A$.
Let $p:A\to P(A)$ be the map given by $p(a)= (a-1)({\sf mod} \
I^3(A))$. Then $p$ is a \emph{quadratic map}, meaning that the
\emph{cross-effect} $$(a\mid b)_p:=p(a+b)-p(a)-p(b)$$ is linear in
$a$ and $b$. Actually $p$ is a universal quadratic function $p:A\to
P(A)$ (see \cite{approximation}). A quadratic map $f:A\to B$ of
abelian groups is called \emph{homogeneous} if $f(-a)=f(a)$. It is
well known that for any abelian group $A$ there  exists a universal
homogeneous quadratic function $\gamma:A\to \Gamma(A)$. If $A$ is a
module over a monoid ring $R=\Z[L]$, then $P(A)$, $\Gamma(A)$, $A\tp
A$  are also $R$-modules, where the action of $x\in L$ is given by
$$p(a)x=p(ax), \ \ (\gamma(a))x=\gamma(ax), \ \ (a\tp b)x=ax\tp bx.$$

\begin{Le}\label{galatoad} If $F\in \bf F$ is an additive functor, then
$$\Hom_{\bf F}(\Gamma\circ \La^2,F)=0=\Hom_{\bf F}(\La^2,F)$$
\end{Le}

{\it Proof}. Let us recall that if $T\in \bf F$ is a functor with $T(0)=0$,
then the second cross-effect of $T$ fits in the decomposition
$$
T(A\oplus B)\cong T(A)\oplus T(B) \oplus T(A\mid B).
$$
Putting $B=A$ and using the codiagonal morphism $(\Id,\Id):A\oplus
A\to A$ one obtains a natural transformation $\eta_A:T(A\mid A)\to
T(A)$. It is clear that any natural transformation from $T$ to an
additive functor factors through ${\sf Coker}(\eta)$. We first take
$T= \Gamma(\La ^2)$. Since the second cross-effect of $\Gamma\circ
\La^2$ contains as a direct summand the term $\Gamma(A\t B)$ and for
$A=B$ it maps via $\eta$ surjectively to $\Gamma(\La^2 A)$, we
conclude that there is no nontrivial map from $\Gamma\circ \La^2$ to
any additive functor. Similarly for $\Hom_{\bf F}(\La^2,F)$.

\begin{Le}\label{pesappr} Let $L$ be a free monoid and let $R=\Z[L]$ be the
corresponding monoid ring. Then
$$\Ext^p_{\bf F}(P,F)=0$$
provided $F$ is additive and $2\le p\le 4$.
\end{Le}
{\it Proof}. Since $R=\Z[L]$ is torsion free as an abelian group and
$F$ is an additive functor the main result of \cite{approximation}
shows that one has an isomorphism
$$\Ext^p_{\bf F}(P,F)\cong \Ext^p_{\bf Q}(P,F), \ \ $$
provided $p\le 4$. Here $\bf Q$ is the abelian category of quadratic
functors from $\modf R$ to $\Mod R$. For the functor $P\tp R$, which
is given by $X\mapsto P(X)\tp R$, one has an isomorphism (see
\cite{approximation})
$$\Hom_{\bf Q}(P\tp R,T)\cong T(R), \ \ T\in \bf Q$$ It follows that $P\tp R$
is a projective object in $\bf Q$. Thus one can use the bar-resolution
$$0\leftarrow P\leftarrow P\tp R\leftarrow P\tp R\tp R\leftarrow \cdots$$
to get a projective resolution of $P$ in the category $\bf Q$. In particular
one has an isomorphism
$$\Ext^*_{\bf Q}(P,F)\cong H^*(R,F(R))$$
and the result follows from the fact that the Hochschild cohomology
vanishes for free rings in  dimensions $\ge 2$.
\subsection{Proof of Theorem \ref{nilnul}} \label{proofnil} We put
 $F=(-)\tp_RB\in \bf
F$. Thanks to Proposition \ref{ssnil} one has the following spectral sequence
$$E^2_{pq}=\Ext^p_{\bf F}(\Tor_q^Q(-,R),F)\Longrightarrow \HML^{p+q}(
Q,B).
$$
By Proposition \ref{qtor} restriction of the functor $\Tor_*^Q(-,R)$
to the category $\modf R$ is given by
$$\Tor_0^Q(-,R)=\Id, $$
$$ \Tor_1^Q(-,R)= \La^2,$$
$$  \Tor_2^Q(-,R)= \Id\t \La^2,$$
$$\Tor_3^Q(-,R)= (\La ^2\t \Id^{\t 2})\oplus (\Gamma\circ \La^2),$$
Since $F$ is additive, Lemma \ref{vanishing} shows that
$$E_2^{p0}= \Ext^p_{{\bf F}}(\Id,F),$$
$$E_2^{p1}= \Ext^p_{{\bf F}}(\La^2,F),$$
$$E_2^{p2}= 0$$
$$E_2^{p3}= \Ext^p_{{\bf  F}}(\Gamma\circ \La^2,F).$$
We also have  $$E^{03}_2=0=E^{01}_2$$ thanks to Lemma \ref{galatoad}. Moreover
 $E_2^{20}=0$ by Lemma \ref{freering}. Thus it suffices to
show that the following differentials of the spectral sequence
$$d_2:E^{11}_2=\Ext^1_{{\bf F}}(\La^2, F)\to E^{30}_2= \Ext^{3}_{{\bf F}}(\Id,F)$$
and
$$d_2:E^{12}_2=\Ext^2_{{\bf F}}(\La^2,F)\to E^{40}_2=\Ext^{4}_{{\bf F}}(\Id,F)$$
are isomorphisms. Let us observe that in general the differential
$$d_2:\Ext^p_{{\bf F}}(\La ^2, F)\to \Ext^{p+2}_{{\bf F}}(\Id,F)$$
is given by the cup product with $e\in \Ext^2_{{\bf F}}(\Id,\La^2)$
corresponding to the extension
$$0\leftarrow \Id\leftarrow  P^2\leftarrow   \Id^{\t 2}\leftarrow \La ^2 \leftarrow  0$$
We have $e=e_1\cup e_2$, where $e_1$ corresponds to the extension
$$0\to \La^2 \to \Id^2\to {\sf Sym}^2\to 0$$
while $e_2$ corresponds to the extension
$$
0\to {\sf Sym}^2\to  P\to \Id\to 0,
$$
where ${\sf Sym}^2$ is the second symmetric power and the first
nontrivial map is induced by the assignment $a\tp b\mapsto (a\mid
b)_p$, while the second map is given by $p(a)\mapsto a.$

It follows from Lemma \ref{vanishing} that the cup product with $e_1$ yields an
isomorphism
$$\Ext^p_{{\bf F}}(\La ^2, F)\to \Ext^{p+1}_{{\bf F}}({\sf Sym ^2},F), \ p\ge -1.$$
 Similarly Lemma \ref{pesappr} shows that the map
$$\Ext^p_{{\bf F}}({\sf Sym} ^2, F)\to \Ext^{p+1}_{{\bf F}}(\Id,F)$$
 induced by the cup product with $e_2$ is an
isomorphism if $2\le p\le 3$ and we are done.

\end{document}